\documentclass{llncs}
\usepackage{eurosym}

\usepackage{times}
\usepackage{amssymb}
\usepackage[toc,page]{appendix}

\usepackage{verbatim}
\usepackage{amsmath}
\usepackage{tikz}
\usepackage{mathrsfs}
\usepackage[margin=10pt,font=small,labelfont=bf]{caption}
\usepackage{dsfont}
\usepackage{enumerate}

\newcounter{cptsx}

\newcommand{\sembrack}[1]{[\![#1]\!]}

\newcommand{\knw}[1]{\textsf{K}}

\renewcommand{\phi}{\varphi}

\usepackage{tikz}
\usetikzlibrary{calc,matrix,positioning,fit, shapes, automata, arrows, patterns}

\tikzset{
  solid node/.style={circle,draw,inner sep=1.2,fill=black},
  empty node/.style={font=\tiny, circle, outer sep=0pt, inner sep=0pt},
  matrix node/.style={rectangle, inner sep=0pt, draw, minimum height=0.5cm, minimum width=0.8cm, outer sep=0pt},  
  noline/.style={edge from parent/.style={draw=none}},
line/.style={edge from parent/.style={draw}},
}

\newcommand{\ben}{\begin{enumerate}\item}
\newcommand{\een}{\end{enumerate}}
\newcommand{\bit}{\begin{itemize}\item}
\newcommand{\eit}{\end{itemize}}

\usepackage{environ}
\makeatletter
\newdimen\royalignsep@
\def\royalign@preamble{%
   &\hfil
    \strut@
    \setboxz@h{\@lign$\m@th\displaystyle{##}$}%
    \ifmeasuring@\savefieldlength@\fi
    \set@field
    \tabskip\z@skip
   &\setboxz@h{\@lign$\m@th\displaystyle{{}##}$}%
    \ifmeasuring@\savefieldlength@\fi
    \set@field
    \hfil
    \tabskip\royalignsep@
}
\NewEnviron{royalign}[1]{%
  \royalignsep@=#1\let\align@preamble=\royalign@preamble
  \begin{align}\BODY\end{align}}
\NewEnviron{royalign*}[1]{%
  \royalignsep@=#1\let\align@preamble=\royalign@preamble
  \begin{align*}\BODY\end{align*}}
\makeatother

\title{A Logic of Objective and Subjective Oughts}

\author{Aldo Iv\'an Ram\'irez Abarca\inst{1}\and Jan Broersen\inst{2}}
\institute{Utrecht University, Utrecht 3512 JK, The Netherlands \email{a.i.ramirezabarca@uu.nl}\and Utrecht University, Utrecht 3512 JK, The Netherlands \email{J.M.broersen@uu.nl}}

\begin{document}

\maketitle

\begin{abstract}

The relation between agentive action, knowledge, and obligation is central to the understanding of responsibility --a main topic in Artificial Intelligence. Based on the view that an appropriate formalization of said relation would contribute to the development of ethical AI, we point out the  main characteristics of a logic for objective and subjective oughts that was recently introduced in the literature. This logic extends the traditional stit paradigm with deontic and epistemic operators, and provides a semantics that deals with Horty's puzzles for knowledge and obligation. We provide an axiomatization for this logic, and address its soundness and completeness with respect to a class of relevant models.


\end{abstract}


\section{Introduction}

AI developers face a big challenge in creating systems that are expected to make ethically charged decisions. The field of machine ethics has seen a quick growth in recent years, and questions regarding \emph{responsibility} of autonomous agents are very important. In our opinion, these questions can be categorized in two main trends: 1) \emph{conceptual} questions that deal with the ontology and essential components of the notion of responsibility, and 2) \emph{technical} questions concerning the implementation of such notion in AI. This work attempts to make a contribution in the \emph{technical} category. We take part in a very specific line of research, where proof systems of deontic logic are intended to help in the testing of ethical behavior of AI through theorem proving (\cite{arkoudas2005toward}, \cite{MUR}, \cite{bringsjord2006toward}). That being said, answers to the technical kind of questions typically presuppose a particular choice of philosophical standpoint against questions of the first kind, the conceptual one. Our philosophical standpoint comes from John Horty's framework of act utilitarian stit logic (\cite{Horty2001}), extended with epistemic relations. As such, we make explicit the goal of the present paper: to provide well-behaved formalizations of 3 essential components of responsibility of intelligent systems: \emph{agentive action}, \emph{knowledge}, and \emph{obligation}. 

According to \cite{arkoudas2005toward}, having an \emph{expressive deontic logic with practical relevance} and an efficient algorithm for proving theorem-hood is highly applicable in the construction of logic-based ethical robots.\footnote{See \cite{pereira2016programming} for an overview of the advantages and disadvantages of doing machine ethics via theorem proving.} To support this statement, Arkoudas et al. present in \cite{arkoudas2005toward} a natural deduction calculus for a logic of ought-to-do that was developed by Horty in \cite{Horty2001} and axiomatised by Murakami in \cite{MUR}. With an interactive theorem proving
system named \textsc{Athena}, they illustrate the fact that we can mechanize deontic logic to do machine ethics.    

The starting point of our work comes from a recent interest in enhancing both the expressivity and practical relevance of Horty's stit theory of ought-to-do in order to deal with situations in which agents' \emph{knowledge} becomes significant. Inspired by 3 puzzles for knowledge-dependent obligations that pose a problem for merely extending his initial logic with epistemic operators, Horty presents in \cite{hortyepistemic} a novel semantics for \emph{epistemic oughts} based on \emph{action types}. Although the approach is substantial, it comes with 3  disadvantages: (1) it diverges from his work in \cite{Horty2001}, (2) there are semantic constraints that limit the expressivity of the models, and (3) the use of action types precludes an efficient axiomatization.\footnote{We will justify all these claims in the second section of this paper, after introducing the 3 puzzles and addressing Horty's solution.} 

With similar motivations as Horty, the extended abstract \cite{JANANDI} proposes an alternative logic --where the main idea is to distinguish objective from subjective obligations-- that allegedly mitigates the disadvantages mentioned above. The authors claim that this logic is simpler, more naturally connected to Horty's work in \cite{Horty2001}, and axiomatizable. Being very brief, \cite{JANANDI} only deals with the conceptual benefits of the proposal, and the proof of soundness and completeness of the logic's proof system is only mentioned. Here, we recover the definition of this logic, address its benefits, and show that its proof system is indeed sound and complete with respect to a class of relevant models. We do this hoping that the results will give some background to new developments in the mechanization of deontic logic for ethical AI, in the aforementioned tradition of \cite{arkoudas2005toward} and \cite{bringsjord2006toward}. 

The paper is structured as follows. After a short presentation of stit and its aplications in the modeling of action, knowledge, and  obligation, we go through Horty's puzzles. We mention his solution to them and justify the claim that his approach with action types comes with disadvantages. Afterwards, we present the logic developed by \cite{JANANDI}, show how it solves Horty's puzzles, and deal with its axiomatizability.\footnote{The proof of completeness is dense and technical, so the full proofs of each statement are provided in the appendix section of the present work.}        

\section{Action, Knowledge, and Obligation in Stit}
\label{shit}

We want to consider oughts from the perspective that an agent should be excused for having failed at an obligation if it lacks the necessary knowledge to perform a required task (\emph{doctors ought to stop the bleeding of the patient, but if they do not know how to, they should be excused}). A typical framework for expressing statements that involve knowledge for required tasks as a component of  responsibility is stit logic (\cite{belnap01facing}). Stit logic was created to formalize the concept of \emph{action}, so it naturally lent itself to the study of \emph{obligation} (\cite{Horty2001}) and of \emph{knowingly doing} (\cite{xstit}), all important elements in the notion of responsibility. For a comprehensive review of the interaction of these 3 concepts in the literature, we refer the reader to \cite{Xu2015}. 

The 3 puzzles of Horty that we mentioned in the introduction are actually very good illustrations of how stit deals with action, knowledge, and obligation, but in order to tackle them we need to recover basic definitions of the paradigm known as \emph{act utilitarian stit logic}. In this paradigm, \emph{obligation} stems from a dominance ordering over the set of available actions. The idea is that the effects of the best actions in the ordering --the so-called \emph{optimal} actions-- are the obligations of a given agent at a given moment. We proceed to introduce the basic aspects of an extension of act utilitarian stit logic with epistemic modalities, and leave its examples to the section where we present Horty's puzzles.   

\begin{definition}[Syntax]
\label{syntax ep stit}
Given a finite set $Ags$ of agent names, a countable set of propositions $P$ such that $p \in P$ and $\alpha \in Ags$, the grammar for the formal language $\mathcal L_{\textsf{KO}}$ is given by:
\[ \begin{array}{lcl}
\phi :=  p \mid \neg \phi \mid \phi \wedge \psi \mid \Box \phi \mid [\alpha] \phi \mid K_\alpha \phi \mid \odot [\alpha] \phi 
\end{array} \]

\end{definition}
\normalsize
$\Box\varphi$ is meant to express the `historical necessity' of $\varphi$
($\Diamond \varphi$ abbreviates $\neg \Box \neg \varphi$). $[\alpha] \varphi$ stands for `agent $\alpha$ sees to it that $\varphi$'. 
$K_\alpha$ is the epistemic operator for $\alpha$. Finally, $\odot [\alpha]\phi$ is meant to represent that $\alpha$ ought to see to it that $\phi$.

As for the semantics, the structures in which we evaluate formulas of the language $\mathcal L_{\textsf{KO}}$  are based on what we call \emph{epistemic act utilitarian branching time frames}. 

\begin{definition}[Branching time (BT) frames]
\label{frames}
A \textbf{finite-choice} epistemic  act utilitarian $BT$-frame is a tuple $\langle T,\sqsubset,\mathbf{\mathbf{Choice}}, \{\sim_\alpha\}_{\alpha\in Ags}, \mathbf{Value}  \rangle$ such
that:

\begin{itemize}

\item $T$ is a non-empty set of \textnormal{moments} and $\sqsubset$ is a strict partial ordering on $T$ satisfying `no backward branching'. Each maximal $\sqsubset$-chain is called a $\textnormal{history}$, which represents a way in which time might evolve. $H$ denotes the set of all histories, and for each $m\in T$, $H_m:=\{h \in H ;m\in h\}$. Tuples $\langle m,h \rangle$ are called \emph{situations} iff $m \in T$, $h \in H$, and $m\in h$.  $\mathbf{Choice}$ is a function that maps each agent $\alpha$ and moment $m$ to a \textbf{finite} partition $\mathbf{Choice}^m_\alpha$ of $H_m$, where the cells of such a partition represent $\alpha$'s available actions at $m$. $\mathbf{Choice}$ satisfies two constraints: \begin{itemize}
\item $(\mathtt{NC})$ or `no choice between undivided histories'.- For all $h, h'\in H_m$, if $m'\in h\cap h'$ for some $m' \sqsupset m$, then $h\in L$ iff $h'\in L$ for every $L\in \mathbf{Choice}^m_\alpha$. 

\item $(\mathtt{IA})$ or `independence of agency'.- A function $s$ on $Ags$ is called a \emph{selection function} at $m$ if it assigns to each $\alpha$ a member of $\mathbf{Choice}^m_\alpha$. If we denote by $\mathbf{Select}^m$ the set of all selection functions at $m$, then we have that for every $m\in T$ and $s\in\mathbf{Select}^m$, $\bigcap_{\alpha \in Ags} s(\alpha)\neq \emptyset$ (see \cite{belnap01facing} for a discussion of the property).
\end{itemize}

\item For $\alpha\in Ags$, $\sim_\alpha$ is the epistemic indistinguishability equivalence relation for agent $\alpha$.


\item $\mathbf{Value}$ is a deontic function that assigns to each history $h\in H$ a real number, representing the utility of $h$.

\end{itemize}
\end{definition}

As mentioned before, the idea is that obligations come from the optimal actions for a given agent. The optimality of such actions is relative to a dominance ordering of the actions, and this ordering is given by the value of the histories in those actions --itself provided by $\mathbf{Value}$. In order to present the semantics for formulas involving the ought-to-do operator, we therefore need some previous definitions. For $\alpha\in Ags$ and $m_*\in T$, we define
a dominance ordering $\preceq$ on $\mathbf{Choice}^ {m_*}_\alpha$ such that for $L, L'\in \mathbf{Choice}^ {m_*}_\alpha$, 
$L\preceq L' \textnormal{ iff }   \mathbf{Value}(h) \leq \mathbf{Value}(h') \textnormal{ for every } h\in L, h'\in L'$. We write $L\prec L'$ \textnormal{iff} $L\preceq L'$ and $L'\npreceq L$. The optimal set of actions, then, is taken as $\mathbf{Optimal}^{m_*}_\alpha:=\{L \in \mathbf{Choice}^{m_*}_\alpha ; \textnormal{there is no } L' \in \mathbf{Choice}^{m_*}_\alpha \textnormal{such that }  L\prec L'\}.$ As is customary, the models and the semantics for the formulas are defined by adding a valuation function to the frames of Definition \ref{frames}: 

\begin{definition}
\label{models KCSTIT}
A BT-model $\mathcal {M}$ consists of the tuple that results from adding a valuation function $\mathcal{V}$ to a BT-frame, where $ \mathcal{V}: P\to 
    2^{T \times H}$ assigns to each atomic proposition a set of moment-history pairs.
    Relative to a model $\mathcal{M}$, the semantics for the formulas of $\mathcal {L}_{\textsf{KO}}$ is defined recursively by the following truth conditions, evaluated at a given situation $\langle m,h \rangle$: 
\[ \begin{array}{lll}
\langle m,h \rangle \models p & \mbox{ iff } & \langle m,h \rangle \in \mathcal{V}(p) \\

\langle m,h \rangle \models \neg \phi & \mbox{ iff } & \langle m,h \rangle \not\models \phi \\

\langle m,h \rangle \models \phi \wedge \psi & \mbox{ iff } & \langle m,h \rangle \models \phi \mbox{ and } \langle m,h \rangle \models \psi \\

\langle m,h \rangle \models \Box \phi &
\mbox{ iff } & \forall h'\in H_m, \langle m,h' \rangle \models \phi \\

\langle m,h \rangle \models [\alpha]
\phi & \mbox{ iff } & \forall h'\in \mathbf{Choice}^m_\alpha(h), \langle m, h'\rangle \models \phi\\

\langle m,h \rangle \models K_{\alpha} \phi &
\mbox{ iff } & \forall \langle m',h'\rangle \mbox{ s.t. }  \langle m,h \rangle \sim_{\alpha}\langle m',h' \rangle, \langle m',h' \rangle \models \phi\\


\langle m,h \rangle \models \odot[\alpha] \phi &
\mbox{ iff } & \forall L\in \mathbf{Optimal} ^{m}_\alpha,  h'\in L \mbox{ implies that }\langle m,h' \rangle \models \varphi.
\end{array} \]
\normalsize
Satisfiability, validity on a frame, and general validity are defined as usual. We write $|\phi|^m$ to refer to the set $\{h\in H_m ;\mathcal{M},\langle m, h\rangle \models \phi\}$.

\end{definition}

\subsection{Horty's Puzzles}
\label{Formal account of the epistemic puzzles}
The 3 puzzles that we have mentioned since the beginning of the paper, and that pose a problem for formalizing epistemic oughts just with the epistemic extension of act utilitarian logic, can be summarized as follows. 

\begin{example}\label{problem1}
Agent $\beta$ places a coin on top of a table --either heads up or tails up-- but hides it from agent $\alpha$. Agent $\alpha$ can bet that the coin is heads up, that it is tails up, or it can refrain from betting. If $\alpha$ bets and chooses correctly, it wins \euro 10. If it chooses incorrectly, it does not win anything, and if it refrains from betting, it wins \euro 5.
\end{example}

\begin{figure}[htb!]
\begin{minipage}[c]{.9\linewidth}
\centering
\resizebox{160 pt}{!}{%
\begin{tikzpicture}[level distance=2cm,
level 1/.style={sibling distance=6cm},
level 2/.style={sibling distance=4.2cm},
level 3/.style={sibling distance=1.5cm}
]

\node {} [grow=up]
	child{node [matrix, matrix of nodes, ampersand replacement=\&, label=left:$m_1$, label=below right:$Choice_{\beta}^{m_1}$] (matrixi) {
	\node[matrix node, label=below:$L_1$] (m11) {}; \& \node[matrix node, label=below:$L_2$] (m12) {}; \\ }  
		child[noline]{node [matrix, matrix of nodes, ampersand replacement=\&, label=right:$m_3$, label=below right:$Choice_\alpha^{m_3}$]{\node[matrix node, label=below:$L_6$] (m31) {};  \& \node[matrix node, label=below:$L_7$] (m32) {}; \& \node[matrix node, label=below:$L_8$] (m33) {};\\}
			child{node (h6) {$h_6$}
				}
			child{node (h5) {$h_5$}
				}
            child{node (h4) {$h_4$}}
			}
		child[noline]{node [matrix, matrix of nodes, ampersand replacement=\&, label=left:$m_2$, label=below left:$Choice_\alpha^{m_2}$]{\node[matrix node, label=below: $L_3$] (m21) {};  \& \node[matrix node, label=below:$L_4$] (m22) {}; \& \node[matrix node, label=below:$L_5$] (m23) {};\\}
			child{node (h3) {$h_3$}
				}
			child{node (h2) {$h_2$}
				}
           	child{node (h1) {$h_1$}
				}
			}
		};

\draw (m11.center) -- (m22);
\draw (m12.center) -- (m32);
\draw (m21.center) -- (h1) node [pos=.7,draw=none, label=above left:\footnotesize$10$] (lh1) {};
\draw (m21.center) -- (h1) node [pos=.5,draw=none, label=above left:\footnotesize$BH$]  {};
\draw (m21.center) -- (h1) node [pos=.3,draw=none, label=above left:\footnotesize$H$] {};
\draw (m22.center) -- (h2) node [pos=.7,draw=none, label=above left:\footnotesize$0$] (lh2) {};
\draw (m22.center) -- (h2) node [pos=.5,draw=none, label=above left:\footnotesize$BT$]  {};
\draw (m22.center) -- (h2) node [pos=.3,draw=none, label=above left:\footnotesize$H$]  {};
\draw (m31.center) -- (h4) node [pos=.7,draw=none, label=above left:\footnotesize$0$] (lh4) {};
\draw (m31.center) -- (h4) node [pos=.5,draw=none, label=above left:\footnotesize$BH$]  {};
\draw (m31.center) -- (h4) node [pos=.3,draw=none, label=above left:\footnotesize$T$]  {};
\draw (m23.center) -- (h3) node [pos=.7,draw=none, label=above left:\footnotesize$5$](lh3) {};
\draw (m23.center) -- (h3) node [pos=.5,draw=none, label=above left:\footnotesize$\lnot G$] {};
\draw (m23.center) -- (h3) node [pos=.3,draw=none, label=above left:\footnotesize$ H$] {};
\draw (m33.center) -- (h6) node [pos=.7,draw=none, label=above left:\footnotesize$5$] (lh6) {};
\draw (m33.center) -- (h6) node [pos=.5,draw=none, label=above left:\footnotesize$\lnot G$] {};
\draw (m33.center) -- (h6) node [pos=.3,draw=none, label=above left:\footnotesize$T$] {};
\draw (m32.center) -- (h5) node [pos=.7,draw=none, label=above left:\footnotesize$10$] (lh5) {};
\draw (m32.center) -- (h5) node [pos=.5,draw=none, label=above left:\footnotesize$BT$] {};
\draw (m32.center) -- (h5) node [pos=.3,draw=none, label=above left:\footnotesize$T$] {};
\draw[dashed,blue, bend left] (m23) to (m31) {};
\end{tikzpicture}
}
\captionof{figure}{Coin problem \#1}
\label{fig1}
\end{minipage}
\end{figure}
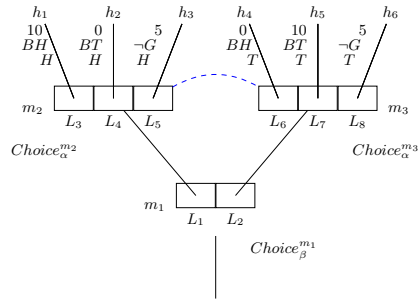

The stit diagram that represents \emph{Horty's interpretation} of the situation is included in Figure \ref{fig1}. In this diagram, we take $H$ to denote the proposition  `$\beta$ places the coin heads up', $T$ to denote `$\beta$ places the coin tails up', $BH$ to denote  `$\alpha$ bets heads', $BT$ to denote  `$\alpha$ bets tails', and $G$ to denote  `$\alpha$ gambles'. In moment $m_1$, $\beta$ places the coin on top of the table, so that its available actions are labeled by $L_1$ (placing the coin heads up) and $L_2$ (placing the coin tails up). At moments $m_2$ and $m_3$, it is $\alpha$'s turn to act, and the available actions are clear from the picture. The blue dotted line represents the epistemic class of $\alpha$: since $\beta$ is hiding the coin, $\alpha$ cannot distinguish whether it is at moment $m_2$ or $m_3$.\footnote{Notice that Horty's formalization also yields that in none of the situations will $\alpha$ knowingly perform any of the available actions: it cannot epistemically distinguish between the situations in which it is `betting heads', `betting tails', or `refraining'.} For such an interpretation regarding the epistemic structure of the agent, a problem ensues due to the fact that for every $i\in \{2,3\}$ and $h\in H_{m_i}$, we have that $\langle m_i,h\rangle \models K_\alpha\odot [\alpha]G$. This means that $\alpha$ knows that it ought to gamble, even if this is a `risky' move that could result in a payoff of $0$. In this sense, we could say that the agent's knowledge of what is optimal would lead it into taking a chance and gambling.

\begin{example}\label{problem2}
With the same scheme as in the previous example, if $\alpha$ bets and chooses correctly, it wins \euro 10. If it refrains from betting, it \emph{also} wins \euro 10. If it bets incorrectly, it does not win anything. 
\end{example}
Intuitively, $\alpha$ ought to refrain from gambling in this scenario, for refraining implies that it would win by the same amount as when betting correctly but without engaging in an action that could possibly fail. In this case, the problem is that for every $i\in \{2,3\}$ and $h\in H_{m_i}$, we have that $\langle m_i,h\rangle \not \models K_\alpha\odot [\alpha]\lnot G$: $\alpha$ does not know that it ought to refrain from gambling. Figure \ref{fig2} includes the stit diagram for this scenario. 

\begin{figure}[htb!]
\begin{minipage}[c]{.9\linewidth}
				\centering
\resizebox{160 pt}{!}{%
\begin{tikzpicture}[level distance=2cm,
level 1/.style={sibling distance=6cm},
level 2/.style={sibling distance=4.2cm},
level 3/.style={sibling distance=1.5cm}
]

\node {} [grow=up]
	child{node [matrix, matrix of nodes, ampersand replacement=\&, label=left:$m_1$, label=below right:$Choice_{\beta}^{m_1}$] (matrixi) {
	\node[matrix node, label=below:$L_1$] (m11) {}; \& \node[matrix node, label=below:$L_2$] (m12) {}; \\ }  
		child[noline]{node [matrix, matrix of nodes, ampersand replacement=\&, label=right:$m_3$, label=below right:$Choice_\alpha^{m_3}$]{\node[matrix node, label=below:$L_6$] (m31) {};  \& \node[matrix node, label=below:$L_7$] (m32) {}; \& \node[matrix node, label=below:$L_8$] (m33) {};\\}
			child{node (h6) {$h_6$}
				}
			child{node (h5) {$h_5$}
				}
            child{node (h4) {$h_4$}}
			}
		child[noline]{node [matrix, matrix of nodes, ampersand replacement=\&, label=left:$m_2$, label=below left:$Choice_\alpha^{m_2}$]{\node[matrix node, label=below: $L_3$] (m21) {};  \& \node[matrix node, label=below:$L_4$] (m22) {}; \& \node[matrix node, label=below:$L_5$] (m23) {};\\}
			child{node (h3) {$h_3$}
				}
			child{node (h2) {$h_2$}
				}
           	child{node (h1) {$h_1$}
				}
			}
		};

\draw (m11.center) -- (m22);
\draw (m12.center) -- (m32);
\draw (m21.center) -- (h1) node [pos=.7,draw=none, label=above left:\footnotesize$10$] (lh1) {};
\draw (m21.center) -- (h1) node [pos=.5,draw=none, label=above left:\footnotesize$BH$]  {};
\draw (m21.center) -- (h1) node [pos=.3,draw=none, label=above left:\footnotesize$H$] {};
\draw (m22.center) -- (h2) node [pos=.7,draw=none, label=above left:\footnotesize$0$] (lh2) {};
\draw (m22.center) -- (h2) node [pos=.5,draw=none, label=above left:\footnotesize$BT$]  {};
\draw (m22.center) -- (h2) node [pos=.3,draw=none, label=above left:\footnotesize$H$]  {};
\draw (m31.center) -- (h4) node [pos=.7,draw=none, label=above left:\footnotesize$0$] (lh4) {};
\draw (m31.center) -- (h4) node [pos=.5,draw=none, label=above left:\footnotesize$BH$]  {};
\draw (m31.center) -- (h4) node [pos=.3,draw=none, label=above left:\footnotesize$T$]  {};
\draw (m23.center) -- (h3) node [pos=.7,draw=none, label=above left:\footnotesize$10$](lh3) {};
\draw (m23.center) -- (h3) node [pos=.5,draw=none, label=above left:\footnotesize$\lnot G$] {};
\draw (m23.center) -- (h3) node [pos=.3,draw=none, label=above left:\footnotesize$ H$] {};
\draw (m33.center) -- (h6) node [pos=.7,draw=none, label=above left:\footnotesize$10$] (lh6) {};
\draw (m33.center) -- (h6) node [pos=.5,draw=none, label=above left:\footnotesize$\lnot G$] {};
\draw (m33.center) -- (h6) node [pos=.3,draw=none, label=above left:\footnotesize$T$] {};
\draw (m32.center) -- (h5) node [pos=.7,draw=none, label=above left:\footnotesize$10$] (lh5) {};
\draw (m32.center) -- (h5) node [pos=.5,draw=none, label=above left:\footnotesize$BT$] {};
\draw (m32.center) -- (h5) node [pos=.3,draw=none, label=above left:\footnotesize$T$] {};
\draw[dashed,blue, bend left] (m23) to (m31) {};
\end{tikzpicture}
}
\captionof{figure}{Coin problem \#2}
\label{fig2}
\end{minipage}
\end{figure}
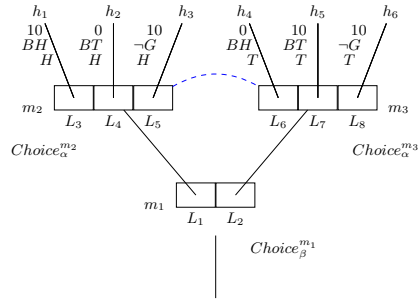
\begin{example}\label{problem3}
With the same scheme as in the previous examples, if $\alpha$ bets and chooses correctly, it wins \euro 10. If it bets incorrectly or refrains from betting, it does not win anything. 
\end{example}

\begin{figure}[htb!]
\begin{minipage}[c]{.9\linewidth}
				\centering
\resizebox{160 pt}{!}{%
\begin{tikzpicture}[level distance=2cm,
level 1/.style={sibling distance=6cm},
level 2/.style={sibling distance=4.2cm},
level 3/.style={sibling distance=1.5cm}
]

\node {} [grow=up]
	child{node [matrix, matrix of nodes, ampersand replacement=\&, label=left:$m_1$, label=below right:$Choice_{\beta}^{m_1}$] (matrixi) {
	\node[matrix node, label=below:$L_1$] (m11) {}; \& \node[matrix node, label=below:$L_2$] (m12) {}; \\ }  
		child[noline]{node [matrix, matrix of nodes, ampersand replacement=\&, label=right:$m_3$, label=below right:$Choice_\alpha^{m_3}$]{\node[matrix node, label=below:$L_6$] (m31) {};  \& \node[matrix node, label=below:$L_7$] (m32) {}; \& \node[matrix node, label=below:$L_8$] (m33) {};\\}
			child{node (h6) {$h_6$}
				}
			child{node (h5) {$h_5$}
				}
            child{node (h4) {$h_4$}}
			}
		child[noline]{node [matrix, matrix of nodes, ampersand replacement=\&, label=left:$m_2$, label=below left:$Choice_\alpha^{m_2}$]{\node[matrix node, label=below: $L_3$] (m21) {};  \& \node[matrix node, label=below:$L_4$] (m22) {}; \& \node[matrix node, label=below:$L_5$] (m23) {};\\}
			child{node (h3) {$h_3$}
				}
			child{node (h2) {$h_2$}
				}
           	child{node (h1) {$h_1$}
				}
			}
		};

\draw (m11.center) -- (m22);
\draw (m12.center) -- (m32);
\draw (m21.center) -- (h1) node [pos=.7,draw=none, label=above left:\footnotesize$10$] (lh1) {};
\draw (m21.center) -- (h1) node [pos=.5,draw=none, label=above left:\footnotesize$W$]  {};
\draw (m21.center) -- (h1) node [pos=.3,draw=none, label=above left:\footnotesize$H$] {};
\draw (m22.center) -- (h2) node [pos=.7,draw=none, label=above left:\footnotesize$0$] (lh2) {};
\draw (m22.center) -- (h2) node [pos=.5,draw=none, label=above left:\footnotesize$BT$]  {};
\draw (m22.center) -- (h2) node [pos=.3,draw=none, label=above left:\footnotesize$H$]  {};
\draw (m31.center) -- (h4) node [pos=.7,draw=none, label=above left:\footnotesize$0$] (lh4) {};
\draw (m31.center) -- (h4) node [pos=.5,draw=none, label=above left:\footnotesize$BH$]  {};
\draw (m31.center) -- (h4) node [pos=.3,draw=none, label=above left:\footnotesize$T$]  {};
\draw (m23.center) -- (h3) node [pos=.7,draw=none, label=above left:\footnotesize$0$](lh3) {};
\draw (m23.center) -- (h3) node [pos=.5,draw=none, label=above left:\footnotesize$\lnot G$] {};
\draw (m23.center) -- (h3) node [pos=.3,draw=none, label=above left:\footnotesize$ H$] {};
\draw (m33.center) -- (h6) node [pos=.7,draw=none, label=above left:\footnotesize$0$] (lh6) {};
\draw (m33.center) -- (h6) node [pos=.5,draw=none, label=above left:\footnotesize$\lnot G$] {};
\draw (m33.center) -- (h6) node [pos=.3,draw=none, label=above left:\footnotesize$T$] {};
\draw (m32.center) -- (h5) node [pos=.7,draw=none, label=above left:\footnotesize$10$] (lh5) {};
\draw (m32.center) -- (h5) node [pos=.5,draw=none, label=above left:\footnotesize$W$] {};
\draw (m32.center) -- (h5) node [pos=.3,draw=none, label=above left:\footnotesize$T$] {};
\draw[dashed,blue, bend left] (m23) to (m31) {};
\end{tikzpicture}
}
\captionof{figure}{Coin problem \#3}
\label{fig3}
\end{minipage}
\end{figure}
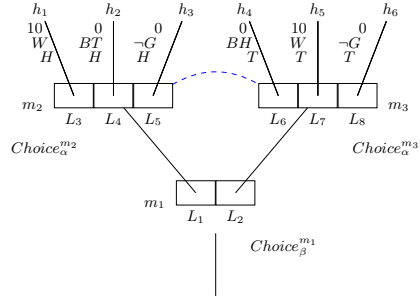

The problem here is that for every $i\in \{2,3\}$ and $h\in H_{m_i}$, $\langle m_i,h\rangle \models K_\alpha  \odot [\alpha]W$, where $W$ is the proposition that stands for `$\alpha$ wins'. This means that $\alpha$ knows that it ought to win at any given situation, but such knowledge is not \emph{action-guiding}, meaning that it will not provide $\alpha$ with a choice to make. Though the agent knowingly ought to win, it cannot knowingly do so --it simply does not have the means due to a lack of knowledge. Thus, Kant's principle of `ought implies can' is not satisfied ($\langle m_i,h\rangle\not\models  K_\alpha  \odot [\alpha]W \to \Diamond  K_\alpha [\alpha]W$). Figure \ref{fig3} includes the stit diagram for this scenario.

Horty solves these 3 puzzles by introducing both syntactic and semantic addenda to epistemic act utilitarian stit logic. He extends the language with an operator $[... \ \mathtt{kstit} ]$ to encode the concept of \emph{ex interim} knowledge (or knowingly doing). The semantics for formulas involving this operator uses \emph{action types}, with the premise that actions of the same \emph{type} may lead to different outcomes in different moments (see \cite{HortyPacuit2017}). Unfortunately, the introduction of the new operator and of types comes with two unfavorable semantic constraints: \begin{enumerate}
    \item  In order for $[... \ \mathtt{kstit} ]$ to be an \textbf{S5} operator, the epistemic relations must ensue not between moment-history pairs but between moments. The problem with this condition is that it limits the class of models to those in which knowledge is moment-dependent (agents will not be able to know that they perform a given action),\footnote{Horty's models satisfy the following constraint: if $\langle m, h\rangle \sim_\alpha \langle m',h'\rangle$ then $\langle m,h'\rangle \sim_\alpha \langle m',h''\rangle$ for every $h'\in H_m$, $h''\in H_{m'}$.}
    \item Indistinguishable moments must offer the same available types. The problem with this constraint is that it cannot be characterized syntactically without producing an infinite axiomatization. This is due to the fact that performing a certain action type can only be expressed syntactically with propositional constants (again, see \cite{HortyPacuit2017} for the details).\footnote{It is an open problem to determine whether there is a finite axiomatization of Horty's logic of epistemic action and obligation if the types were also included in the object language.}
\end{enumerate}  

As a solution to the stated puzzles, Horty's approach is successful. However, \cite{JANANDI} claims that we can also be successful without using action types. The benefits of the framework presented in \cite{JANANDI}, then, include technical characterizability of the  constraints imposed on the structures --which is important for axiomatization--, semantic simplicity, and enhanced expressivity. In the following section, all these claims will be substantiated. 

\section{A Logic of Objective and Subjective oughts}

\cite{JANANDI} proposes to disambiguate two senses of ought-to-do in order to produce a system that solves Horty's puzzles, avoids action types, and is axiomatizable. The two senses are an \emph{objective} one, which coincides with Horty's act utilitarian ought-to-do, and a \emph{subjective} one, which arises from the epistemically best candidates in the set of available actions for a given agent. By `epistemically best' we mean those actions that are undominated not only in the actual moment but whose all epistemic equivalents across different indistinguishable situations are also undominated. 

Essentially, we are talking about an extension of the language in  Definition \ref{syntax ep stit} with a new operator $\odot_{\mathcal{S}}[\alpha]$, meant to build up formulas that would express what $\alpha$ subjectively ought to do. As for the semantics of this new operator, it involves a dominance ordering as well, but one different to that which is used for objective ought-to-do's. In order to define this subjective dominance ordering, \cite{JANANDI} introduces a new semantic concept known as \emph{epistemic clusters}, which are nothing more than a given action's epistemic equivalents in situations that are indistinguishable to the actual one. Formally, we have that for  $\alpha\in Ags$, $m_*, m\in T$, and $L\subseteq H_{m_*}$, $L$'s \emph{epistemic cluster} at $m$ is the set \[[L]^m_\alpha:=\{h\in H_m ; \exists h_*\in L \ \textnormal{ s.t. }\ \langle m_*,h_*\rangle \sim_\alpha \langle m,h\rangle \}.\] 
As a convention, we write $m \sim_\alpha m'$ if there exist $h\in H_m$, $h'\in H_{m'}$ such that $\langle m,h\rangle \sim_\alpha \langle m',h'\rangle$. The notion of epistemic clusters is used to define a subjective dominance ordering $\preceq_s$ on $\mathbf{Choice}^{m_*}_\alpha$ by the following rule: for $L, L'\in\mathbf{Choice}^{m_*}_\alpha$, $L\preceq_s L'$ iff for every $m$ such that $m_*\sim_\alpha m$, $\mathbf{Value}(h) \leq \mathbf{Value}(h')$ for every $h\in [L]^m_\alpha, h'\in [L']^m_\alpha$. Just as in the case of objective ought-to-do's, this ordering allows us to define a subjectively optimal set of actions $\mathbf{S-optimal}^{m_*}_\alpha:=\{L \in \mathbf{Choice}^{m_*}_\alpha ; \textnormal{ there is no } L' \in \mathbf{Choice}^{m_*}_\alpha \textnormal{ s. t. }  L\prec_s L'\},$ where we write $L\prec_s L'$ iff $L\preceq_s L'$ and $L'\npreceq_s L$. The idea, then, is that something will be a subjective obligation of a given agent at a given moment if it is an effect of all the subjectively optimal actions of that agent at that moment.   

As established in \cite{JANANDI}, the models in which to evaluate the formulas of the extended language need to satisfy extra constraints in order to capture an appropriate interaction of action, knowledge, and subjective obligation. By this we mean to say that in these models (a) agents should be able to knowingly do the same things across epistemically indistinguishable states, (b) the subjective ought-to-do must conform to Kant's directive of `ought implies can' in its epistemic version of `subjectively ought-to-do implies ability of knowingly doing', and (c) if something is a subjective ought-to-do of a given agent, then the agent should know that that is the case. Therefore, we focus on models that fulfill the following requirements, which will grant the conditions mentioned above (as can be seen from the proof of soundness): \begin{itemize}
\item $(\mathtt{OAC})$ For every situation $\langle m_*,h_*\rangle$, if $\langle m_*, h_*\rangle\sim_\alpha \langle m, h\rangle$ for some $\langle m,h\rangle$, then $\langle m_*,h_*'\rangle\sim_\alpha \langle m,h\rangle$ for every $h_*'\in \mathbf{Choice}^{m_*}_\alpha (h_*)$. We refer to this constraint as the `own action condition', since it implies that agents do not know more than what they perform. Because of this constraint, the knowledge that we are formalizing here is of a very particular kind: to know something is just the same as to knowingly do it.   

\item $(\mathtt{Unif-H})$ For every situation $\langle m_*,h_*\rangle$, if $\langle m_*, h_*\rangle \sim_\alpha \langle m, h\rangle $ for some $\langle m,h\rangle $, then for every $h_*'\in H_{m_*}$, there exists $h'\in H_m$ such that $\langle m_*, h_*'\rangle \sim_\alpha \langle m, h'\rangle$. Combined with $(\mathtt{OAC})$, this constraint is meant to capture a notion of uniformity of strategies, where epistemically indistinguishable situations should offer similar actions for the agent to choose upon. We call this condition `uniformity of historical possibility'. 
\end{itemize}
In finite-choice epistemic act utilitarian BT-models that satisfy these two constraints, then, the semantics for the formulas involving $\odot_{\mathcal{S}}[\alpha]$ is defined as expected: \[\langle m,h \rangle \models \odot_{\mathcal{S}}[\alpha] \varphi \mbox{ iff }  \forall L\in \mathbf{S-optimal} ^{m}_\alpha, \forall m' \mbox{ s.t. } m\sim_\alpha m', \ [L]^{m'}_\alpha\subseteq |\varphi|^{m'}.\]
\subsection{Solution to Horty's Puzzles}
The semantics for subjective ought-to-do's offers solutions to natural interpretations of Horty's puzzles, in which the assumption that the coin is hidden from the betting agent is captured by taking $\sim_\alpha$ to be defined by the following information sets: $\{\langle m_2, h_1\rangle, \langle m_3, h_4\rangle\}$, in which $\alpha$ bets heads; $\{\langle m_2, h_2\rangle, \langle m_3, h_5\rangle\}$, in which $\alpha$ bets tails; and $\{\langle m_2, h_3\rangle , \langle m_3, h_6\rangle\}$, in which $\alpha$ refrains from betting. 

\begin{figure}[hbt!]
    \begin{minipage}[c]{.9\linewidth}
    \centering
\resizebox{160 pt}{!}{%
\begin{tikzpicture}[level distance=2cm,
level 1/.style={sibling distance=6cm},
level 2/.style={sibling distance=4.2cm},
level 3/.style={sibling distance=1.5cm}
]

\node {} [grow=up]
	child{node [matrix, matrix of nodes, ampersand replacement=\&, label=left:$m_1$, label=below right:$Choice_{\beta}^{m_1}$] (matrixi) {
	\node[matrix node, label=below:$L_1$] (m11) {}; \& \node[matrix node, label=below:$L_2$] (m12) {}; \\ }  
		child[noline]{node [matrix, matrix of nodes, ampersand replacement=\&, label=right:$m_3$, label=below right:$Choice_\alpha^{m_3}$]{\node[matrix node, label=below:$L_6$] (m31) {};  \& \node[matrix node, label=below:$L_7$] (m32) {}; \& \node[matrix node, label=below:$L_8$] (m33) {};\\}
			child{node (h6) {$h_6$}
				}
			child{node (h5) {$h_5$}
				}
            child{node (h4) {$h_4$}}
			}
		child[noline]{node [matrix, matrix of nodes, ampersand replacement=\&, label=left:$m_2$, label=below left:$Choice_\alpha^{m_2}$]{\node[matrix node, label=below: $L_3$] (m21) {};  \& \node[matrix node, label=below:$L_4$] (m22) {}; \& \node[matrix node, label=below:$L_5$] (m23) {};\\}
			child{node (h3) {$h_3$}
				}
			child{node (h2) {$h_2$}
				}
           	child{node (h1) {$h_1$}
				}
			}
		};

\draw (m11.center) -- (m22);
\draw (m12.center) -- (m32);
\draw (m21.center) -- (h1) node [pos=.7,draw=none, label=above left:\footnotesize$10$] (lh1) {};
\draw (m21.center) -- (h1) node [pos=.5,draw=none, label=above left:\footnotesize$BH$]  {};
\draw (m21.center) -- (h1) node [pos=.3,draw=none, label=above left:\footnotesize$H$] {};
\draw (m22.center) -- (h2) node [pos=.7,draw=none, label=above left:\footnotesize$0$] (lh2) {};
\draw (m22.center) -- (h2) node [pos=.5,draw=none, label=above left:\footnotesize$BT$]  {};
\draw (m22.center) -- (h2) node [pos=.3,draw=none, label=above left:\footnotesize$H$]  {};
\draw (m31.center) -- (h4) node [pos=.7,draw=none, label=above left:\footnotesize$0$] (lh4) {};
\draw (m31.center) -- (h4) node [pos=.5,draw=none, label=above left:\footnotesize$BH$]  {};
\draw (m31.center) -- (h4) node [pos=.3,draw=none, label=above left:\footnotesize$T$]  {};
\draw (m23.center) -- (h3) node [pos=.7,draw=none, label=above left:\footnotesize$5$](lh3) {};
\draw (m23.center) -- (h3) node [pos=.5,draw=none, label=above left:\footnotesize$\lnot G$] {};
\draw (m23.center) -- (h3) node [pos=.3,draw=none, label=above left:\footnotesize$ H$] {};
\draw (m33.center) -- (h6) node [pos=.7,draw=none, label=above left:\footnotesize$5$] (lh6) {};
\draw (m33.center) -- (h6) node [pos=.5,draw=none, label=above left:\footnotesize$\lnot G$] {};
\draw (m33.center) -- (h6) node [pos=.3,draw=none, label=above left:\footnotesize$T$] {};
\draw (m32.center) -- (h5) node [pos=.7,draw=none, label=above left:\footnotesize$10$] (lh5) {};
\draw (m32.center) -- (h5) node [pos=.5,draw=none, label=above left:\footnotesize$BT$] {};
\draw (m32.center) -- (h5) node [pos=.3,draw=none, label=above left:\footnotesize$T$] {};
\draw[dashed,blue, bend right] (lh1) to (lh4) {};
\draw[dashed,blue, bend right] (lh2) to (lh5) {};
\draw[dashed,blue, bend right] (lh3) to (lh6) {};

\end{tikzpicture}
}
\captionof{figure}{Coin problem \#1, revisited}
\label{figura}
\end{minipage}
\end{figure}
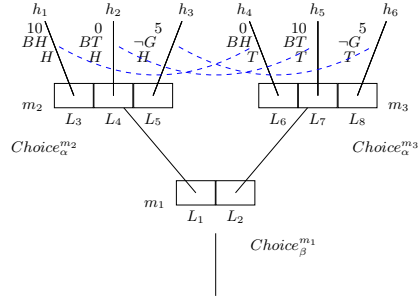

For Example \ref{problem1}, the problem is solved because although $\langle m_i,h_i\rangle \models K_\alpha\odot [\alpha]G$, we consider this as the knowledge of an \emph{objective} ought-to-do. \emph{Subjectively} speaking, we do not obtain that $\alpha$ knows that it ought to gamble: $\langle m_i, h_i\rangle \not\models \odot_{\mathcal{S}}[\alpha] G$ and thus $\langle m_i, h_i\rangle \not\models K_\alpha \odot_{\mathcal{S}}[\alpha] G$. This can be seen by noticing that $\mathbf{S-Optimal}^{m_2}_\alpha=\{L_3, L_4, L_5\}$ and $\mathbf{S-Optimal}^{m_3}_\alpha=\{L_6, L_7, L_8\}$. Figure \ref{figura} provides a stit diagram for this scenario.

In Example \ref{problem2}, the problem is solved because $\langle m_i,h_i\rangle\models \odot_{\mathcal{S}}[\alpha] \lnot G$ and $\langle m_i,h_i\rangle\models K_\alpha\odot_{\mathcal{S}}[\alpha] \lnot G$  (notice that $\mathbf{S-Optimal}^{m_2}_\alpha=\{L_5\}$ and $\mathbf{S-Optimal}^{m_3}_\alpha=\{L_8\}$). 

For Example \ref{problem3}, the problem is solved because we obtain that $\langle m_i,h_i\rangle \not\models \odot_{\mathcal{S}}[\alpha]  W$, which in turn implies that $\langle m_i,h_i\rangle \not\models K_\alpha\odot_{\mathcal{S}}[\alpha]  W$  (notice that $\mathbf{S-Optimal}^{m_2}_\alpha=\{L_3, L_4, L_5\}$ and $\mathbf{S-Optimal}^{m_3}_\alpha=\{L_6, L_7, L_8\}).$ Therefore, although $\alpha$ knows that it objectively ought to win, it is not the case that it subjectively ought to win.

When comparing these solutions to Horty's, it is important to point out that the formalization we are using is different from his, for it presupposes that the indistinguishability relation for agent $\alpha$ ensues between \emph{situations}. Regardless of such difference, the solution is virtually the same: for Example \ref{problem1}, Horty gets that for every $h\in H_{m_2}$ $\langle m_2, h\rangle  \not\models \odot[\alpha \ \mathtt{kstit}] G$. For Example \ref{problem2}, he gets that for every $h\in H_{m_2}$ $\langle m_2, h\rangle  \models \odot[\alpha \ \mathtt{kstit}] \lnot G$, and for Example \ref{problem3}, he gets that for every $h\in H_{m_2}$ $\langle m_2, h\rangle  \not\models \odot[\alpha \ \mathtt{kstit}] W$. Therefore, we can see that the notion of $\odot_{\mathcal{S}}[\alpha]$ works as an analog of $\odot [... \ \mathtt{kstit} ]$.

\section{Axiomatization and Some Logical Properties}

In this section, we introduce a proof system for the logic presented, address its soundness and completeness results, and mention a few interesting properties of it.

\begin{definition}[Proof system]
\label{axiomsystemurakami}
Let $\Lambda$ be the proof system defined by the following axioms and rules of inference:

\noindent \textit{(Axioms)} \begin{itemize}
\item All classical tautologies from propositional logic.

\item The $\mathbf{S5}$ axiom schemata for $\square$, $[\alpha ]$, $K_\alpha$,\footnote{The $\mathbf{S5}$ axiom schemata are standard from modal logic, and we include their names here just for coherence. They are $(K)$, $(T)$, $(4)$, and $(5)$.}
\item The following axiom schemata for the interactions of formulas with the given operators:
\small
\begin{flalign*}
&\odot [\alpha ] (p\to q)\to (\odot [\alpha ] p \to \odot [\alpha ] q)& (A1)\\ 
&\square p\to [\alpha ] p \land \odot [\alpha ] p& (A2)\\ 
&\square\odot [\alpha ] p\lor  \square\lnot\odot [\alpha ] p&(A3)\\ 
& \odot [\alpha ] p\to  \odot [\alpha ] ([\alpha ]p)&(A4)\\
&\odot [\alpha ] p\to \Diamond  [\alpha ] p &(Oic)\\
&\mbox{For $n\geq 1$ and pairwise different $\alpha_1,\dots,\alpha_n$},& \\ &\bigwedge_{1\leq k\leq n}  \Diamond [\alpha_i ] p_i \to \Diamond\left(\bigwedge_{1\leq k\leq n}[\alpha_i ] p_i\right)& (IA)\\
&\odot_{\mathcal{S}}[\alpha ] (p\to q)\to (\odot_{\mathcal{S}} [\alpha ] p \to \odot_{\mathcal{S}}[\alpha ] q)& (A5)\\
&\odot_{\mathcal{S}} [\alpha ] p\to  \odot_{\mathcal{S}}[\alpha ] (K_\alpha p)&(A6)\\
&K_\alpha p\to [\alpha ]p&(OAC)\\
&\Diamond K_\alpha p \to K_\alpha \Diamond  p&(Unif-H)\\
&K_\alpha \square p\to\odot_{\mathcal{S}}[\alpha ]p&(s.N)\\
&\odot_{\mathcal{S}}[\alpha ]p\to\Diamond K_\alpha p&(s.Oic)\\
&\odot_{\mathcal{S}}[\alpha ] p \to K_{\alpha} \square\odot_{\mathcal{S}}[\alpha ] p  &(Cl)
\end{flalign*}

\end{itemize} 
\noindent \textit{(Rules of inference)} 
\begin{itemize} \item \textit{Modus Ponens}, Substitution, and Necessitation for the modal operators.
\end{itemize}
\end{definition}

Schema $(IA)$ encodes `independence of agency'. $(OAC)$ encodes the `own action condition'. $(Unif-H)$ encodes the `uniformity of historical possibility' constraint. $(Oic)$ and $(s.Oic)$ concern the objective, resp. subjective, versions of Kant's directive of `ought implies can'. $(s.N)$ (standing for `subjective necessity') captures that all that is historically necessary at epistemically indistinguishable situations must be a subjective obligation, and $(Cl)$ (standing for `closure') characterizes a property that says that if one subjectively ought to do something, then one knows that that is settled. 

The axiom system $\Lambda$ turns out to be sound and complete with respect to the class of epistemic act utilitarian \emph{bi-valued} BT-models. Such models are more general than the ones introduced in Definition \ref{models KCSTIT}. Instead of only one value function, there are two: one for the objective ought-to-do's, and the other for the subjective ones.\footnote{This extension is extremely useful for the proof of completeness, but neither its conceptual reach nor its philosophical implications has been a subject of our investigation as of yet. There are more than a few reasons to entertain skepticism about this technical decision, since it is not grounded philosophically. To have two deontic functions instead of one lends itself to potentially having different values for histories according to whether their utility is objective or subjective. However, we must observe that the notions of objective and subjective utility are not what drives the disambiguation of ought-to-do into objective and subjective obligations. What drives it is the epistemic structure of a given agent. Therefore, the use of two value functions  should not be taken as essential for the proposal of the logic presented here. Moreover, the soundness result works for models in which there is only one deontic function, and the solution to Horty's puzzles does not rely on there being two. On a related note, we make explicit that the focus of this work was to develop well-behaved semantics for subjective ought-to-do's in the stit tradition, so that if we presented a language with only the subjective-ought operator and not the objective one, as well as with the operators for knowledge, action, and historical necessity, then our proof of completeness would in fact be enough to show completeness of this fragment with respect to the class of single-valued epistemic act utilitarian BT-models.} These models, then, are of the form $\langle T, \sqsubset, \mathbf{Choice},\{\sim_\alpha\}_{\alpha\in Ags}, \mathbf{Value}_{\mathcal{O}}, \mathbf{Value}_{\mathcal{S}} \rangle$. As such, the models in Definition \ref{models KCSTIT} are particular instances of bi-valued models, in which both value functions assign the same value to each history of the tree. 

Furthermore, our soundness and completeness results presuppose a logic of ought-to-do that deals with the `sure-thing principle', according to which the ranking of the available actions of a given agent should take into consideration what all the other agents choose concurrently. For a given agent, the action profiles of the other agents are seen as the possible states in which the agent will act (see \cite{Horty2001}, Chapter 4, subsection 4.1.2). We implement Horty's approach for dealing with the  `sure-thing principle' in the logic that we axiomatize, and this means that the semantics that we use is in fact a generalization of the one introduced before. In order to address such a generalization and the soundness and completeness results, we need further definitions:

For $m\in T$ and $\beta\in Ags$, we define \[\mathbf{State}_\beta^m=\{S\subseteq H_m; S=\bigcap_{\alpha \in Ags-\{\beta\}} s(\alpha), \mbox{ where }s\in \mathbf{Select}^m\}.\] With this notion of states, we redefine the dominance orderings so that the actions are measured taking into consideration the states which those actions are facing. For $\alpha\in Ags$ and $m_*\in T$, we first define
a general ordering $\leq$ on $\mathcal{P}(H_{m_*})$ such that for $X, Y\subseteq H_{m_*}$, $
X\leq Y \textnormal{ iff } \mathbf{Value}_\mathcal{O}(h) \leq \mathbf{Value}_\mathcal{O}(h') \textnormal{ for every } h\in X, h'\in Y$. The objective dominance ordering $\preceq$ is now defined such that for $L, L'\in \mathbf{Choice}_\alpha^{m_*}$, $
L\preceq L' \textnormal{ iff }  \forall S\in \mathbf{State}_\alpha^{m_*}, L\cap S \leq L'\cap S.$ 
As for the subjective ought-to-do's, we first define a general ordering $\leq_s$ on $\mathcal{P}(H_{m_*})$ such that for $X, Y\subseteq H_{m_*}$, 
$X\leq_s Y \textnormal{ iff }\mathbf{Value}_\mathcal{S}(h) \leq_s \mathbf{Value}_\mathcal{S}(h') \textnormal{ for all } h\in X, h'\in Y.$ The subjective dominance ordering $\preceq_s$ is then defined  such that for $L, L'\subseteq H_{m_*}$, $L\preceq_s L'$ iff $\forall m$ such that $m_*\sim_\alpha m$, $\forall S\in \mathbf{State}_\alpha^{m}, [L]^m_\alpha\cap S \leq_s [L']^m_\alpha\cap S.$

With these definitions, we adapt the semantics for the formulas involving ought-to-do operators, making the logic strong enough to deal with the `sure-thing principle': we set $\langle m,h \rangle \models \odot[\alpha] \varphi$ iff for every $L\in \mathbf{Choice}^{m}_\alpha$ such that $\langle m,h_L\rangle\nvDash\varphi$ for some $h_L\in L$, there exists $L'\in \mathbf{Choice}_\alpha^{m}$ such that $L\prec L'$ and if $L''=L$ or $L'\preceq L''$, then $\langle m,h'\rangle\vDash \varphi$ for every $h'\in L''$. Similarly, we set $\langle m,h \rangle \models \odot_{\mathcal{S}}[\alpha] \varphi$ iff for every $L\in \mathbf{Choice}^{m}_\alpha$ such that $\langle m',h_L\rangle\nvDash\varphi$ for some $h_L\in [L]^{m'}_\alpha (m\sim_\alpha m')$, there exists $L'\in \mathbf{Choice}_\alpha^{m}$ such that $L\prec_s L'$ and if $ L''=L'$ or $L'\preceq_s L''$, then for every $m''$  such that $m\sim_\alpha m'', \langle m'',h''\rangle\vDash \varphi$ for every $h''\in [L'']^{m''}_\alpha.$

\subsection{Soundness}
\begin{proposition}\label{soundness} The system $\Lambda$ is sound with respect to the class of epistemic act utilitarian bi-valued BT-models. 
\end{proposition}

\begin{proof}
Standard. See Appendix for details. 
\end {proof}

\subsection{Completeness}

For reasons of space, we can only summarize the proof here, leaving the full detailed exposition for the Appendix.  The proof of completeness of $\Lambda$ with respect to the class of epistemic act utilitarian bi-valued BT-models is a two-step process. First, we introduce relational structures called Kripke-estit models for evaluating the formulas of the language $\mathcal L_{\textsf{KO}}$ under a semantics that mirrors the semantics of Definition \ref{models KCSTIT}, and prove completeness of $\Lambda$ with respect to these structures via the technique of canonical models. Secondly, we provide a truth-preserving correspondence between Kripke-estit models and certain epistemic act utilitarian bi-valued BT-models, so that completeness with respect to Kripke-estit models grants it with respect to bi-valued BT-models. Essential to the success of the technique of canonical models is to prove the so-called \emph{existence} and \emph{truth} lemmas just as in modal logic. In the case of this particular proof system, these lemmas are quite involved, and the reader is encouraged to go over them in the Appendix, which offers enough room to address all their features carefully and with precision. In our opinion, the results of soundness and completeness are significant, all the more because we want to provide some theoretical foundations to the idea --explored by \cite{arkoudas2005toward}-- that we can actually do machine ethics via theorem-proving (or model-checking). 

To end on a slightly less technical note, we finish this section by addressing some interesting properties concerning interactions of certain operators in the logic presented:
\begin{equation*}
\not\models  \odot[\alpha] \varphi \to \Diamond K_{\alpha} [\alpha]\varphi. 
\end{equation*}

\emph{`It is not the case that if the agent objectively ought to do something, then it can knowingly do it}.' Our solution to Example 3 poses a counterexample, because $\langle m_2, h_1\rangle$ is such that $\langle m_2, h_1\rangle \models\odot [\alpha]W $, but $\langle m_2, h_1\rangle \not\models\Diamond K_\alpha [\alpha] W$, as witnessed by the facts that $\langle m_2,h_2\rangle \not\models [\alpha]W$ and that $\langle m_3,h_4\rangle \not\models [\alpha ]W$.
\begin{equation*}\not\models  \odot_{\mathcal{S}}[\alpha] \varphi \to  \odot[\alpha] \varphi.  \end{equation*}
\emph{`It is not the case that if the agent subjectively ought to do something, then it objectively ought to do it}.'
Our solution to Example 2 poses a counterexample,  because $\langle m_2, h_1\rangle$ is such that $\langle m_2, h_1\rangle\models\odot_{\mathcal{S}} [\alpha] \lnot G $, but $\langle m_2, h_1\rangle\not\models \odot [\alpha] \lnot G$.
\begin{equation*}\not\models \odot[\alpha] \varphi \to  \odot_{\mathcal{S}}[\alpha] \varphi. \end{equation*}
\emph{`It is not the case that if the agent objectively ought to do something, then it subjectively ought to do it}.'
Our solution to Example 3 poses a counterexample, because $\langle m_2, h_1\rangle$ is such that $\langle m_2, h_1\rangle\models\odot [\alpha] W $, but $\langle m_2, h_1\rangle\not\models\odot_{\mathcal{S}} [\alpha] W$.

\section{Conclusion}

This work deals with important questions in the modeling of agency, knowledge, and obligation. Formal depictions of such concepts are likely to be useful when it comes to doing machine ethics based on deontic logic and its mechanization. The approach discussed and analyzed here is based on a stit logic of utilitarian `ought-to-do' enriched with epistemic relations. We argue that to solve certain problems in the treatment of knowledge and obligations within stit --namely Horty's puzzles-- one possibility is to distinguish between objective and subjective versions of the ought-to-do modality. Moreover, we show that this possibility comes with formal advantages such as simplicity of semantics and axiomatizability.

\bibliography{AAAI2018references}

\begin{thebibliography}{10}
\providecommand{\url}[1]{\texttt{#1}}
\providecommand{\urlprefix}{URL }
\providecommand{\doi}[1]{https://doi.org/#1}

\bibitem{arkoudas2005toward}
Arkoudas, K., Bringsjord, S., Bello, P.: Toward ethical robots via mechanized
  deontic logic. In: AAAI Fall Symposium on Machine Ethics. pp. 17--23 (2005)

\bibitem{belnap01facing}
Belnap, N., Perloff, M., Xu, M.: Facing the future: agents and choices in our
  indeterminist world. Oxford University Press (2001)

\bibitem{bringsjord2006toward}
Bringsjord, S., Arkoudas, K., Bello, P.: Toward a general logicist methodology
  for engineering ethically correct robots. IEEE Intelligent Systems
  \textbf{21}(4),  38--44 (2006)

\bibitem{xstit}
Broersen, J.: Deontic epistemic stit logic distinguishing modes of mens rea.
  Journal of Applied Logic  \textbf{9}(2),  137--152 (2011)

\bibitem{JANANDI}
Broersen, J., Ram{\'\i}rez~Abarca, A.I.: Formalising oughts and practical
  knowledge without resorting to action types. In: Proceedings of the 17th
  International Conference on Autonomous Agents and MultiAgent Systems. pp.
  1877--1879. International Foundation for Autonomous Agents and Multiagent
  Systems (2018)

\bibitem{hortyepistemic}
Horty, J.: Epistemic oughts in stit semantics  (2018)

\bibitem{HortyPacuit2017}
Horty, J., Pacuit, E.: Action types in stit semantics. Review of Symbolic Logic
   (2017), forthcoming

\bibitem{Horty2001}
Horty, J.F.: Agency and Deontic Logic. Oxford University Press (2001)

\bibitem{MUR}
Murakami, Y.: Utilitarian deontic logic. AiML-2004: Advances in Modal Logic
  \textbf{287} (2004)

\bibitem{pereira2016programming}
Pereira, L.M., Saptawijaya, A., et~al.: Programming machine ethics, vol.~26.
  Springer (2016)

\bibitem{Xu2015}
Xu, M.: Combinations of stit with ought and know. Journal of Philosophical
  Logic  \textbf{44}(6),  851--877 (2015). \doi{10.1007/s10992-015-9365-7},
  \url{http://dx.doi.org/10.1007/s10992-015-9365-7}

\end{thebibliography}
\bibliographystyle{splncs04}
\clearpage

\begin{appendix}
\section{Proof of Soundness}
\begin{proposition} The system $\Lambda$ is sound with respect to the class of epistemic utilitarian bi-valued BT-models. 
\end{proposition}

\begin{proof}
Let $\mathcal{M}=\langle T,\sqsubset,\mathbf{\mathbf{Choice}}, \{\sim_\alpha\}_{\alpha\in Ags}, \mathbf{Value}_\mathcal{O}, \mathbf{Value}_\mathcal{S}  \rangle$ be an epistemic utilitarian bi-valued BT-frame.  Let $\mathcal{V}$ be any pertinent valuation function. The $\mathbf{S5}$ axioms for $\square, [\alpha ]$, $K_\alpha$, as well as axioms $(A1)$, $(A2)$, $(A3)$, $(A4)$, $(A5)$, $(A6)$, $(Oic)$, $(IA)$, $(OAC)$, and $(Unif-H)$ are shown to be semantic validities straightforwardly. Since they involve some novelty, we include the detailed proofs for $(s.N)$, $(s.Oic)$, and $(Cl)$ here.

\begin{itemize}
\item To see that $\mathcal{M}\vDash (s.N)$, take $\langle m_*,h_*\rangle$ such that $\langle m_*,h_*\rangle\models K_\alpha \square \varphi$. Let $m$ be such that $m_*\sim_\alpha m$ (which means that there exist $j\in H_{m_*}$, $j'\in H_{m}$ such that $\langle m_*,j\rangle \sim_\alpha \langle m,j'\rangle$), and let $L\in \mathbf{Choice}^{m_*}_\alpha$. Condition ($\mathtt{Unif-H}$) ensures that there exists $h_*'\in H_m$ such that $\langle m_*,h_*\rangle \sim_\alpha \langle m,h_*'\rangle$. Our assumption that  $\langle m_*,h_*\rangle\models K_\alpha \square \varphi$ implies then that $\langle m,h_*'\rangle\models \square \varphi$. Therefore, for any $h\in [L]^m_\alpha$, the fact that $h\in H_m$ yields that $\langle m,h\rangle\models \varphi$. Therefore, for every $L\in \mathbf{Choice}^{m_*}_\alpha$ and every $m$ such that $m_*\sim_\alpha m$ , we have that $[L]^m_\alpha\subseteq |\varphi|^m$, which vacuously implies that $\langle m_*,h_*\rangle\models \odot_{\mathcal{S}}[\alpha ] \varphi$.

\item To see that $\mathcal{M}\vDash (s.Oic)$, take $\langle m_*,h_*\rangle$ such that $\langle m_*,h_*\rangle\models \odot_{\mathcal{S}}[\alpha ] \varphi$. This implies that there exists $L_*\subseteq H_{m_*}$ such that $[L_*]^{m'}_\alpha \subseteq |\varphi|^{m'}$ for every $m'\in T$ such that $m_*\sim_\alpha m'$. Since $\sim_\alpha$ is reflexive, it is clear that $[L_*]^{m_*}_\alpha \subseteq |\varphi|^{m_*}$. Now, take $h_{0}\in L_*$. Let $\langle m, h \rangle $ be such that $\langle m_*, h_{0}\rangle \sim_\alpha \langle m,h\rangle$. From the definition of the epistemic cluster set, we have that $h\in [L_*]^{m}_\alpha$, so that $[L_*]^{m}_\alpha \subseteq |\varphi|^{m}$ and thus $\langle m,h\rangle \models \varphi$.  Therefore, we have that   $h_{0}\in H_{m_*}$ is such that for every $\langle m, h\rangle$ with $\langle m_*, h_{0}\rangle \sim_\alpha \langle m,h\rangle$, $\langle m,h\rangle\models \varphi$. This means that $\langle m_*,h_{0}\rangle \models K_\alpha \varphi$. Therefore, we have that $\langle m_*,h_{*}\rangle \models \Diamond K_\alpha \varphi$.

\item To see that $\mathcal{M}\vDash(Cl)$, take $\langle m_*,h_*\rangle$  such that $\langle m_*,h_*\rangle\models \odot_{\mathcal{S}}[\alpha ] \varphi$. Let $\langle m,j\rangle$ be such that $\langle m_*,h_*\rangle\sim_\alpha \langle m,j\rangle$. Take $h\in H_m$. We want to show that for every $L\in \mathbf{Choice}^{m}_\alpha$ such that $[L]^{m'}\not\subseteq|\varphi|^{m'}$ (for some $m'$ such that $m\sim_\alpha m'$), there exists $L'\in \mathbf{Choice}_\alpha^{m} \mbox{ such that } L\prec_s L' \mbox{ and if } L''=L'\mbox{ or } L'\preceq_s L'', \mbox{ then for every } m'' \mbox{ such that } m\sim_\alpha m''$, we have that $[L'']^{m''}_\alpha\subseteq |\varphi|^{m''}.$  Take $L\in \mathbf{Choice}^{m}_\alpha$ such that there exists $m'\in T$ with $m\sim_\alpha m'$ and $[L]^{m'}\not\subseteq|\varphi|^{m'}$. Let $N_L$ be an action in $\mathbf{Choice}^{m_*}_\alpha$ such that $N_L\subseteq [L]^{m_*}_\alpha$. We know that such an action exists, in virtue of $(\mathtt{Unif-H})$ and $(\mathtt{OAC})$. Notice that transitivity of $\sim_\alpha$ entails that $[N_L]^{o}_\alpha = [L]^{o}_\alpha$ for any moment $o$, so we have that $[N_L]^{m'}_\alpha\not\subseteq|\varphi|^{m'}$. Since  $\langle m_*,h_*\rangle\models \odot_{\mathcal{S}}[\alpha ] \varphi$, we get that there must exist $N\in \mathbf{Choice}^{m_*}_\alpha$ such that $N_L\prec_s N \mbox{ and if } N'=N\mbox{ or } N\preceq_s N', \mbox{ then for every } m'' \mbox{ such that } m_*\sim_\alpha m'', [N']^{m''}_\alpha\subseteq |\varphi|^{m''}.$ Now, let $L_{N}$ be an action in $\mathbf{Choice}^{m}_\alpha$ such that $L_N\subseteq [N]^{m}_\alpha$ (which implies that $[L_N]^{o}_\alpha = [N]^{o}_\alpha$ for any moment $o$). We claim that $L\prec_s L_N$, and show our claim with the following argument. Let $m''$ be a moment such that $m\sim_\alpha m''$, and let $S\in \mathbf{State}^{m''}_\alpha$. On one hand, we have that $\star$ $[L]^{m''}_\alpha\cap S =[N_L]^{m''}_\alpha\cap S\leq_s [N]^{m''}_\alpha\cap S=[L_N]^{m''}_\alpha\cap S$. On the other hand, we know that there exists a moment $m'''$ such that $m_*\sim_\alpha m'''$ and a state $S_0\in \mathbf{State}^{m'''}_\alpha$ such that $[N]^{m'''}_\alpha\cap S_0\not\leq_s[N_L]^{m'''}_\alpha\cap S_0$. Therefore, we have that  $(\star)$  $[L_N]^{m'''}_\alpha\cap S_0=[N]^{m'''}_\alpha\cap S_0\not\leq_s[N_L]^{m'''}_\alpha\cap S_0=[L]^{m'''}_\alpha\cap S_0$. Together, $(\star)$ and $(\star\star)$ entail that $L\prec_s L_N$. Now, let $L''\in \mathbf{Choice}^m_\alpha$ such that $L''=L_N$ or $L_N\preceq_sL''$. If $L'' =L_N$, then for every $m''$ such that $m\sim_\alpha m''$, $[L'']^{m''}_\alpha=[N]^{m''}_\alpha\subseteq |\varphi|^{m''}$. If $L_N\prec_sL''$, then an argument similar to the the one we used above to show that our claim was true renders that an action $N_{L''}\in \mathbf{Choice}^{m_*}_\alpha$ with $N_{L''}\subseteq [L'']^{m_*}_\alpha$ is such that $N\preceq_s N_{L''}$, so that $[L'']^{m''}_\alpha=[N_{L''}]^{m''}_\alpha\subseteq |\varphi|^{m''}$. With this, we have shown that $\langle m,h\rangle\models \odot_{\mathcal{S}}[\alpha ]\varphi$. Since this happens for every $h\in H_m$, we have that $\langle m,j\rangle\models \square\odot_{\mathcal{S}}[\alpha ] \varphi$. Since this happens for any $\langle m,j\rangle$ such that $\langle m_*,h_*\rangle\sim_\alpha \langle m,j\rangle$, then we get that $\langle m_*,h_*\rangle\models K_\alpha\square\odot_{\mathcal{S}}[\alpha ] \varphi$.

\item \emph{Modus Ponens}, Substitution, and Necessitation for $\square$ preserve validity. 
\end{itemize}
Therefore, we have shown that the system $\Lambda$ is sound with respect to the class of epistemic bi-valued BT-models.
\end {proof}

\section{Proof of Completeness}

\begin{definition}[Kripke-estit frames]\label{ledbetter}

A tuple $\mathcal{F}=\langle W, R_\square, \mathtt{Choice}, \{\mathtt{\approx}_\alpha\}_{\alpha\in Ags}, \mathtt{Value}_\mathcal{O}, \mathtt{Value}_\mathcal{S} 
\rangle$ is called a Kripke-estit frame iff $W$ is a set of possible worlds and

\begin{enumerate}

\item $R_\square$ is an equivalence relation over $W$. For $w\in W$, the class of $w$ under $R_\square$ is denoted by $\overline{w}$.

\item $\mathtt{Choice}$ is a function that assigns to each $\alpha\in Ags$ and each $\square$-class $\overline{w}$ (with $w\in W$) a partition $\mathtt{Choice}_\alpha^{\overline{w}}$ of $\overline{w}$ given by an equivalence relation which we will denote by $R_\alpha^{\overline{w}}$. $\mathtt{Choice}$ must satisfy the following constraint:
\begin{itemize}
\item $\mathtt{(IA)_K}$ For $w\in W$, we have that each function $s:Ags\to \mathcal{P}(\overline{w})$ that maps $\alpha$ to a member of $\mathtt{Choice}^{\overline{w}}_\alpha$ is such that $\bigcap_{\alpha \in Ags} s(\alpha) \neq \emptyset$ (the set of all functions $s$ that map $\alpha$ to a member of $\mathtt{Choice}^{\overline{w}}_\alpha$ is denoted by $\mathtt{Select}^m$). 

\end{itemize}
For $w\in W$ and $v\in \overline{w}$, the class of $v$ in the partition $\mathtt{Choice}^{\overline{w}}_\alpha$ is denoted by $\mathtt{Choice}^{\overline{w}}_\alpha(v)$. For $w\in W$ and $\beta\in Ags$, we define $\mathtt{State}_\beta^{\overline{w}}=\{S\subseteq \overline{w} ;S=\bigcap_{\alpha \in Ags-\{\beta\}} s(\alpha), \mbox{ where } s\in \mathtt{Select}^m\}$.

\item For each $\alpha \in Ags$, $\approx_\alpha$ is an (epistemic) equivalence relation on $W$ that satisfies the following constraints: 
\begin{itemize}
\item $\mathtt{(OAC)_K}$ For $w\in W$ such that $v\in \overline{w}$, if $v\approx_\alpha u$, then $v'\approx_\alpha u$ for every $v'\in \mathtt{Choice}^{\overline{w}}_\alpha(v)$.
\item $\mathtt{(Unif-H)_K}$ Let $w_1, w_2\in W$ such that there exist $v\in \overline{w_1}$ and $u\in \overline{w_2}$ with $v\approx_\alpha u$. Then for every $v'\in \overline{w_1}$, there exists $u'\in \overline{w_2}$ such that $v'\approx_\alpha u'$.
\end{itemize}

For each $v\in W$ such that $w\sim_\alpha v$, we define the epistemic cluster set of $L$ at $\overline{v}$ by $\sembrack{L}_\alpha^{\overline{v}}:=\{u\in {\overline{v}} ;\exists o\in L \ \textnormal{ such that }\ o \approx_\alpha u \}.$

\item $\mathtt{Value}_\mathcal{O}$ and $\mathtt{Value}_\mathcal{S}$ are (utilitarian) deontic functions that assign to each $w\in W$ a value in $\mathds{R}$. We use them to define objective and subjective orderings just as for BT-models. 
\end{enumerate}
\end{definition}

Kripke-estit frames allow us to evaluate the formulas of $\mathcal L_{\textsf{KO}}$ in an \emph{analogous} fashion to the semantics provided for BT-frames. Therefore, we use the name `Kripke-estit models' for the structures that result from adding a valuation function to Kripke-estit frames, where the semantics for the formulas of $\mathcal L_{\textsf{KO}}$ is defined so as to mirror Definition \ref{models KCSTIT}. Kripke-estit models can be used for constructing epistemic utilitarian bi-valued BT-models such that both satisfy the same formulas of $\mathcal{L}_{\textsf{KO}}$:

\begin{definition}[Associated BT-frame] \label{tree}

Let $\mathcal{F}=\langle W, R_\square, \mathtt{Choice},  \{\mathtt{\approx}_\alpha\}_{\alpha\in Ags}, \mathtt{Value}_\mathcal{O}, \mathtt{Value}_\mathcal{S}
\rangle$ be a Kripke-estit frame. We define the tuple $\langle T_W, \sqsubset, \mathbf{Choice},\{\sim_\alpha\}_{\alpha\in Ags}, \mathbf{Value}_{\mathcal{O}}, \mathbf{Value}_{\mathcal{S}} \rangle$ will be called the epistemic utilitarian bi-valued BT-frame associated to $\mathcal{F}$ such that
\begin{enumerate}
\item $T_W:={W}\cup \{\overline{w} ;w\in W \} \cup \{W\}$ and $\sqsubset$ is a relation on $T_W$ such that $\sqsubset$ is defined as the transitive closure of the union $\{(\overline{w},v) ; w\in W \mbox{ and } v\in\overline{w} \}\cup \{(W,\overline{w}) ; w\in W\}$. Notice that $\sqsubset$ is a strict partial order on $T_W$ that satisfies ``\emph{no backward branching}'' straightforwardly. Since the tuple $\langle T_W, \sqsubset\rangle$ is thus a tree, we can call the maximal $\sqsubset$-chains in $T_W$ \emph{histories}. We observe that the definition of $\sqsubset$ yields that there is a bijective  correspondence between the set of worlds in $W$ and the histories of $T_W$. If for each $v\in W$, we set $h_v$ to be the history $\{W, \overline{v},v\}$, then we have that for any world $o\in W$, $o\in h_v$ iff $o=v$. Therefore, we will conveniently identify these histories using the worlds at their terminal nodes. In this line of thought, we have that for each $w\in W$, the set $H_{\overline{w}}$ of histories running through $\overline{w}$ is equal to the set $\{h_v ;v\in \overline{w}\}$ --since $\overline{w}\in h_v$ iff $v\in \overline{w}$. Similarly, we have that the set $H$ of all histories is the same as $\{h_v ;v\in W\}$.

\item $\mathbf{Choice}$ is an agentive choice function for the tree of the above item. In order to define such a function, we use the partitions $\mathtt{Choice}_\alpha^{\overline{w}}$ (where $w\in W$ and $\alpha\in Ags$) in the following manner: for each $B\in \mathcal{P}(W)$, we will write $B^T$ to denote the set $\{h_v ;\ v\in B\}$; this way, we define a function $\mathbf{Choice}$ on $Ags\times T_W$ with extension given by $\mathbf{Choice}(\alpha,v)=\{\{h_v\}\}$ for each $\alpha\in Ags$ and $v\in W$, $\mathbf{Choice}(\alpha, \overline{w})=\{C_\alpha^T ;C_\alpha\in \mathtt{Choice}_\alpha^{\overline{w}}\}$ for each $w\in W$ and each $\alpha\in Ags$, and $\mathbf{Choice}(\alpha, W)=\{H\}$ for each $\alpha\in Ags$. To keep notation consistent, we will denote sets of the form $\mathbf{Choice}(\alpha, \overline{w})$ by $\mathbf{Choice}_\alpha^{\overline{w}}$, so that the choice-class of a given $h_v$ in $\mathbf{Choice}_\alpha^{\overline{w}}$ is denoted by $\mathbf{Choice}_\alpha^{\overline{w}}(h_v)$, as expected.\footnote{\label{gruta} Notice that this implies that, for two worlds $v, v'$ in $\overline{w}$, $vR_\alpha^{\overline{w}} v'$ iff $h_v\in \mathbf{Choice}_\alpha^{\overline{w}}(h_{v'})$.} It is important to mention that this definition yields that for every $S\in \mathtt{State}^{\overline{w}}_\alpha$, $S^T\in \mathbf{State}^{\overline{w}}_\alpha$, and that for each $U\in \mathbf{State}^{\overline{w}}_\alpha$, there exists $V\in\mathtt{State}^{\overline{w}}_\alpha$ such that $U=V^T$.

\item 
For each $\alpha \in Ags$, we define an epistemic relation $\thicksim_\alpha$ on $T_W$ using the set $\{\approx_\alpha\}_{\alpha\in Ags}$. For a given agent $\alpha\in Ags$, we define $\sim_\alpha=\{(\langle\overline{w},h_v\rangle, \langle\overline{w'},h_{v'}\rangle) ;w, w'\in W \mbox{ and } v\approx_\alpha v'\} \cup \{(\langle z,h_z\rangle,\langle z,h_z\rangle) ;z\in W\}\cup\{(\langle W,h_v\rangle,\langle W,h_{v'}\rangle) ;v,v'\in W\}$. We observe, then, that this definition entails that $\sim_\alpha$ is an equivalence relation for every $\alpha\in Ags$ and that, for a given $w\in W$ and $L\in\mathtt{Choice}_\alpha^{\overline{w}}$, $v\in \sembrack{L}_\alpha^{\overline{w}}$ iff $h_v\in [L^T]^{\overline{w}}_\alpha$.

\item $\mathbf{Value}_{\mathcal{O}}$ and $\mathbf{Value}_{\mathcal{S}}$ are two utilitarian deontic functions, defined such that for each $h_v\in H$, $\mathbf{Value}_{\mathcal{O}}(h_v)=\mathtt{Value}_{\mathcal{O}}(v)$ and  $\mathbf{Value}_{\mathcal{S}}(h_v)=\mathtt{Value}_{\mathcal{S}}(v)$. Endowed with these deontic functions, we can define objective and subjective dominance orderings for this structure accordingly.
\end{enumerate}

\end{definition} 

\begin{proposition}\label{avion}
Let $\mathcal{F}=\langle W, R_\square, \mathtt{Choice}, \{\mathtt{\approx}_\alpha\}_{\alpha\in Ags}, \mathtt{Value}_\mathcal{O}, \mathtt{Value}_\mathcal{S} 
\rangle$ be a Kripke-estit frame. The tuple $\langle T_W, \sqsubset, \mathbf{Choice},\{\sim_\alpha\}_{\alpha\in Ags}, \mathbf{Value}_{\mathcal{O}}, \mathbf{Value}_{\mathcal{S}} \rangle$ is an epistemic utilitarian bi-valued BT-frame. 
\end{proposition}

\begin{proof}
It amounts to showing that the tuple, as defined, validates the semantic conditions $(\mathtt{NC})$, $(\mathtt{IA})$, $(\mathtt{OAC})$, and $(\mathtt{Unif-H})$.\begin{itemize}
\item $(\mathtt{NC})$ is vacuously validated in moment $W$. It is validated in moments of the form $\overline{w}$ ($w\in W$), since two different histories never intersect in a moment later than $\overline{w}$). Finally, it is also validated in moments of the form $v$ such that $v\in W$ (since there are no moments above $v$).
\item For $(\mathtt{IA})$, we reason by cases: \begin{enumerate}[a)]
\item At moment $W$, $(\mathtt{IA})$ is validated straightforwardly, since $\mathbf{Choice}(\alpha, W)=\{H\}$ for each $\alpha\in Ags$.
\item For a fixed moment of the form $\overline{w}$ (with $w\in W$), let $s$ be a function that assigns to each agent $\alpha$ a member of $\mathbf{Choice}^{\overline{w}}_\alpha=\{(C_\alpha)^T; C_\alpha\in \mathtt{Choice}^{\overline{w}}_\alpha\}$. Let $s_k:Ags \to \bigcup_{\alpha\in Ags}\mathtt{Choice}^{\overline{w}}_\alpha$ be a function such that $s_k(\alpha)=C_\alpha$ iff $s(\alpha)=(C_\alpha)^T$. Since $\mathcal{M} $ validates the constraint $(\mathtt{IA})_K$, then we have that $\bigcap_{\alpha \in Ags} s_k(\alpha) \neq \emptyset$. Take $v\in \bigcap_{\alpha \in Ags} s_k(\alpha)$. Then $v\in C_\alpha$ for every $\alpha\in Ags$. This implies that $h_v\in (C_\alpha)^T$ for every $\alpha\in Ags$, so that $\bigcap_{\alpha \in Ags} s(\alpha) \neq \emptyset$.
\item At moments of the form $v$ such that $v\in W$, if $s$ is a function that assigns to each agent $\alpha$ a member of $\mathbf{Choice}(v,\alpha)$, we have that $s$ must be constant and $\bigcap_{\alpha \in Ags} s(\alpha) = \{h_v\}$.
\end{enumerate}
\item For $(\mathtt{OAC})$, again we reason by cases: 

\begin{enumerate}[a)]
\item Assume that $\langle \overline{w}, h_v \rangle \sim_\alpha \langle \overline{w'}, h_{v'} \rangle$ (for $w, w' \in W$). This means that $v\approx_\alpha v'$. We want to show that for every $h_u\in \mathbf{Choice}^{\overline{w}}_\alpha$ such that $h_u\in\mathbf{Choice}^{\overline{w}}_\alpha(h_v)$, $\langle \overline{w}, h_u \rangle \sim_\alpha \langle \overline{w'}, h_{v'} \rangle$. Therefore, let $h_u\in\mathbf{Choice}^{\overline{w}}_\alpha(h_v)$. By definition, we have that this means that $u\in\mathtt{Choice}^{\overline{w}}_\alpha(v)$. Since $\mathcal{M} $ validates the constraint $(\mathtt{OAC})_K$, this last fact implies, with $v\approx_\alpha v'$, that $u\approx_\alpha v'$, which in turn yields that $\langle \overline{w}, h_u \rangle \sim_\alpha \langle \overline{w'}, h_v \rangle$.

\item For situations with moments of the form $v$ such that $v\in W$, we have that $(\mathtt{OAC})$  is met straightforwardly, since for $h_v$, the choice-cell in $\mathbf{Choice}(\alpha,v)$ to which $h_v$ belongs is just $\{h_v\}$. 

\item For situations with moment $W$, $(\mathtt{OAC})$ is also met straightforwardly, since for every $\alpha\in Ags$, $\sim_\alpha$ is defined such that $\langle W,h_v\rangle\sim_\alpha \langle W, h_{v'}\rangle$ for every pair of histories $h_v, h_{v'}$ in $H$. 

\end{enumerate}
\item  For $(\mathtt{Unif-H})$, again we reason by cases: \begin{enumerate}[a)]
\item Assume that $\langle \overline{w}, h_v \rangle \sim_\alpha \langle \overline{w'}, h_{v'} \rangle$ (for $w, w' \in W$). This means that $v\in\overline{w}$, $v'\in\overline{w'}$, and $v\approx_\alpha v'$. Let $h_{z}\in H_{\overline{w}}$ (which means that $z\in \overline{w}$). We want to show that there exists $h\in H_{\overline{w'}}$ such that $\langle \overline{w}, h_z \rangle \sim_\alpha \langle \overline{w'}, h \rangle$. Condition $(Unif-H)_K$ gives us that there exists $z'\in \overline{w'}$ such that $z\approx_\alpha z'$, which by definiton of $\sim_\alpha$ means that $\langle \overline{w}, h_z \rangle \sim_\alpha \langle \overline{w'}, h_{z'} \rangle$.

\item For situations with moments of the form $v$ such that $v\in W$, we have that $\langle v, h_v \rangle$ is $\sim_\alpha$-related only to itself, so that $(\mathtt{Unif-H})$  is met straightforwardly. 

\item For situations with moment $W$, $(\mathtt{Unif-H})$ is also met straightforwardly, since $\langle W,h_v\rangle\sim _\alpha\langle W,h_{v'}\rangle$ for every  $v,v'\in W$.
\end{enumerate}
\end{itemize}

\end{proof}

Let $\mathcal{M} $ be a Kripke-estit model with valuation function $\mathcal{V}$. The frame upon which $\mathcal{M} $ is based has an associated BT-frame. If to the tuple of this BT-frame we add a valuation function $\mathcal{V}^T$ such that $\mathcal{V}^T(p)=\{\langle\overline{w},h_w\rangle; w\in \mathcal{V}(p)\}$, we call the resulting model the epistemic utilitarian bi-valued BT-model associated to $\mathcal{M}$.  

\begin{proposition}\label{avion1}
Let $\mathcal{M}$ be a Kripke-estit model, and let $\mathcal{M}^t$ denote its associated BT-model. For a formula $\phi$ of $\mathcal{L}_{\textsf{KO}}$ and $w\in W$, we have that $\mathcal{M} ,w\Vdash\phi \mbox{ iff } \mathcal{M} ^T,\langle\overline{w},h_w\rangle\vDash\phi$.  
\end{proposition}

\begin{proof}

\textbf{Claim 1.}
For a given Kripke-estit model $\mathcal{M} $ and a fixed $w\in W^\Lambda$, we have that for every $L, N\in \mathtt{Choice}^{\overline{w}}_\alpha$,
\begin{itemize}
\item $L\preceq N$ iff $L^T\preceq N^T$ and $L\prec N$ iff  $L^T\prec N^T$.
\item $L\preceq_s N$ iff $L^T\preceq_s N^T$ and $L\prec_s N$ iff   $L^T\prec_s N^T$.
\item $L\in \mathtt{Optimal}^{\overline{w}}_\alpha$ iff $L^T\in \mathbf{Optimal}^{\overline{w}}_\alpha$, and $L\in \mathtt{S-Optimal}^{\overline{w}}_\alpha$ iff $L^T\in \mathbf{S-Optimal}^{\overline{w}}_\alpha$.
\end{itemize}
We will prove \textbf{Claim 1.}only for the strict orderings, since these proofs include the arguments needed to show that the statements for $\preceq$ and $\preceq_s$ also hold. 
\begin{itemize}
\item $(\Rightarrow)$ We assume that $L\prec N$. Let $U\in \mathbf{State}^{\overline{w}}_\alpha$. We know that $U=V^T$ for some $V\in \mathtt{State}^{\overline{w}}_\alpha$. Therefore, we have that $L\cap V\leq N\cap V$. But notice that, by definition, we have that for each $v\in \overline{w}$, $\mathbf{Value}_{\mathcal{O}}(h_v)=\mathtt{Value}_{\mathcal{O}}(v)$, so the fact that $L\cap V\leq N\cap V$ implies that $L^T\cap U\leq N^T\cap U$. Now, our assumption also yields that there should exist $S_0\in \mathtt{State}^{\overline{w}}_\alpha$ such that $N\cap S_0\nleq L\cap S_0$. This straightforwardly implies that $N^T\cap S_0^T\nleq L^T\cap S_0^T$, so that indeed we have that $L^T\prec N^T$.  

$(\Leftarrow)$ We assume that $L^T\prec N^T$. Let $S\in \mathtt{State}^{\overline{w}}_\alpha$. We know that 
$L^T\cap S^T\leq N^T\cap S^T$, so that $L\cap S\leq N\cap S$. Similarly, here our assumption yields that there should exist $U_0\in \mathbf{State}^{\overline{w}}_\alpha$ such that $N^T\cap U_0\nleq L^T\cap U_0$. But we have that $U_0=V_0^T$ for some $V_0\in \mathtt{State}^{\overline{w}}_\alpha$, so that the definition of $\mathbf{Value}_{\mathcal{O}}$ implies that $N\cap V_0\nleq L\cap V_0$. With this, we have shown that $L\prec N$.

\item $(\Rightarrow)$ We assume that $L\prec_s N$. Let $w'$ be such that $w\approx_\alpha w'$, and let $U\in \mathbf{State}^{\overline{w'}}_\alpha$. We know that $U=V^T$ for some $V\in \mathtt{State}^{\overline{w'}}_\alpha$. Therefore, we have that $\sembrack{L}^{\overline{w'}}_\alpha\cap V\leq_s \sembrack{N}^{\overline{w'}}_\alpha\cap V$. We have that for each $w\in W$, $\mathbf{Value}_{\mathcal{S}}(h_w)=\mathtt{Value}_{\mathcal{S}}(w)$, and that for each pair $v_1, v_2$ such that $v_1\in \sembrack{L}^{\overline{w'}}_\alpha\cap V$ and $v_2\in \sembrack{N}^{\overline{w'}}_\alpha\cap V$, then $h_{v_1}\in [L^T]^{\overline{w'}}_\alpha\cap V^T$ and $h_{v_2}\in [N^T]^{\overline{w'}}_\alpha\cap V^T$, so that the fact that $\sembrack{L}^{\overline{w'}}_\alpha\cap V\leq_s \sembrack{N}^{\overline{w'}}_\alpha\cap V$ then implies that $[L^T]^{\overline{w'}}_\alpha\cap U\leq_s [N]^{\overline{w'}}_\alpha\cap U$. Now, our assumption also yields that there should exist $w_*$ such that $w\approx_\alpha w_*$ and $S_0\in \mathtt{State}^{\overline{w_*}}_\alpha$ such that $\sembrack{N}^{\overline{w_*}}_\alpha\cap S_0^T\nleq_s \sembrack{L}^{\overline{w_*}}_\alpha\cap S_0^T$. Analogously, this implies that $[N^T]^{\overline{w_*}}_\alpha\cap S_0\nleq_s [L^T]^{\overline{w_*}}_\alpha\cap S_0$, so that indeed we have that $L^T\prec_s N^T$.  

$(\Leftarrow)$ We assume that $L^T\prec_s N^T$. Let $w'$ be such that $w\approx_\alpha w'$, and let $S\in \mathtt{State}^{\overline{w'}}_\alpha$. We know that 
$[L^T]^{\overline{w'}}_\alpha\cap S^T\leq_s [N^T]^{\overline{w'}}_\alpha\cap S^T$, so that an argument similar to the one displayed for the above item renders that $\sembrack{L}^{\overline{w'}}_\alpha\cap S\leq_s \sembrack{N}^{\overline{w'}}_\alpha\cap S$. On the other hand, there should exist $w_*$ such that $w\approx_\alpha w_*$ and $U_0\in \mathbf{State}^{\overline{w_*}}_\alpha$ such that $N^T\cap U_0\nleq_s L^T\cap U_0$. But we have that $U_0=V_0^T$ for some $V_0\in \mathtt{State}^{\overline{w_*}}_\alpha$, so that, in an analogous fashion to what we did above, we get that $N\cap V_0\nleq_s L\cap V_0$. With this, we have shown that $L\prec_s N$.
\item Straightforward from the above items. 
\end{itemize}

With that we conclude the proof of \textbf{Claim 1.} For the proof of the main statement, we proceed by induction on $\phi$. For the base case, we have that, for a propositional letter $p$ and an arbitrary $w\in W$, $\mathcal{M} ,w\Vdash p$ iff $w\in \mathcal{V}(p)$ iff $\langle\overline{w},h_w\rangle\in \mathcal{V}^T(p)$ iff $\mathcal{M} ^T,\langle\overline{w},h_w\rangle\vDash p$. The cases with the boolean connectives are standard, so let us deal with the modal operators. Let $w\in W$ and $\alpha\in Ags$.

\begin{itemize}

\item (``$\square$'')

We have that $\mathcal{M} ,w\Vdash \square\phi$ iff for every $v\in \overline{w}$, $\mathcal{M}_{w},v\Vdash \phi$, which by induction hypothesis happens iff $\mathcal{M} ^T,\langle\overline{v},h_v\rangle\vDash \phi$ for every $v\in \overline{w}$, which happens iff $\mathcal{M} ^T,\langle\overline{w},h_w\rangle\vDash \square\phi$, since we have that $h_v\in H_{\overline{w}}$ iff $v\in \overline{w}$. 

\item (``$[\alpha ]$'')

We have that $\mathcal{M} ,w\Vdash [\alpha ]\phi$ iff for every $v\in W$ such that $wR_\alpha^{\overline{w}} v$, $\mathcal{M} ,v\Vdash \phi$, which by induction hypothesis happens iff $\mathcal{M} ^T,\langle\overline{w},h_v\rangle\vDash \phi$ for every $h_v\in\mathbf{Choice}^{\overline{w}}_\alpha (h_w)$, which in turn  happens iff $\mathcal{M} ^T,\langle\overline{w},h_w\rangle\vDash [\alpha ]\phi$.

\item (``$K_\alpha$'')

We have that $\mathcal{M} ,w\Vdash K_\alpha \phi$ iff for every $v\in W$ such that $w\approx_\alpha v$, $\mathcal{M} ,v\Vdash \phi$, which by induction hypothesis occurs iff $\mathcal{M} ^T,\langle\overline{v},h_v\rangle\vDash \phi$ for every $h_v\in H$ such that $\langle\overline{w},h_w\rangle\sim_\alpha \langle\overline{v},h_v\rangle$, which happens iff $\mathcal{M} ^T,\langle\overline{w},h_w\rangle\vDash K_\alpha\phi$.

\item (``$\odot_\mathcal[\alpha ]$'')

($\Rightarrow$) We assume that $\mathcal{M} ,w\Vdash \odot[\alpha ]\phi$. Then for every $L \in \mathtt{Choice}_\alpha^{\overline{w}}$ such that $\mathcal{M} ,v_L\nVdash \phi$ for some $v_L\in L$,
there exists $L'\in \mathtt{Choice}_\alpha^{\overline{w}} \mbox{ such that } L\prec L' \mbox{ and if } L''=L'\mbox{ or } L'\preceq L'', \mbox{ then } \mathcal{M} ,v\Vdash \phi \mbox{ for every } v\in L''$. Let $N\in \mathbf{Choice}_\alpha^{\overline{w}}$ such that $\mathcal{M}^T,\langle\overline{w},h_z\rangle\nvDash \phi \mbox{ for some } h_z\in N$. By induction hypothesis, we get that $\mathcal{M} ,z\nVdash\phi$. The cell $N$ is $L_0^T$ for some $L_0\in\mathtt{Choice}^{\overline{w}}_\alpha$. Since $z\in L_0$ and $\mathcal{M} ,z\nVdash\phi$, our assumption entails that there exists $L_1\in \mathtt{Choice}_\alpha^{\overline{w}}$ such that $L_0\prec L_1$ and if $L''=L_1$ or $L_1\preceq L''$ then $\mathcal{M} ,v\Vdash \phi$ for every $v\in L''$. We claim that $L_1^T$ is the choice cell at moment $\overline{w}$ that will witness to the fact that $\mathcal{M} ^T,\langle\overline{w},h_w\rangle\vDash \odot[\alpha ] \phi$. In order to show this, we first notice that the first item of \textbf{Claim 1} renders that $L_0^T\prec L_1^T$. Now, we also have that for every $h_v\in L_1^T$, $v\in L_1$ (observe also that $v\in \overline{w}$ and therefore $\mathcal{M} ,v\Vdash \phi$; in this way, we have that for every $h_v\in L_1^T$, the induction hypothesis gives that $\mathcal{M} ^T,\langle\overline{w},h_v\rangle\vDash \phi$. Finally, let $V^T\in \mathbf{Choice}^{\overline{w}}_\alpha$ such that $L_1^T\preceq V^T$. By the first item of \textbf{Claim 1}, this last condition implies that $L_1\preceq V$. Take $h_v\in V^T$. We know that $v\in V$, so that the fact that $L_1\preceq V$ implies with our assumption that $\mathcal{M} ,v\Vdash \phi$; therefore,  the induction hypothesis gives that $\mathcal{M} ^T,\langle\overline{w},h_v\rangle\vDash \phi$. With this, we have shown that our claim is true and that $\mathcal{M} ^T,\langle\overline{w},h_w\rangle\vDash \odot [\alpha ] \phi$.       

($\Leftarrow$) We assume that $\mathcal{M} ^T,\langle\overline{w},h_w\rangle\vDash \odot[\alpha ] \phi$. Then for every $L^T \in \mathbf{Choice}_\alpha^{\overline{w}}$ such that $\mathcal{M} ^T,\langle\overline{w},h_z\rangle\nvDash \phi$ for some $h_z\in L^T$,
there exists $L^T_1\in \mathbf{Choice}_\alpha^{\overline{w}}$ such that $L^T\prec ^T_1$ and if  $U=L_1^T$ or  $L^T\preceq U$, then $\mathcal{M} ^T,\langle\overline{w},h_u\rangle\vDash \phi \mbox{ for every } h_u\in U$. Let $N\in \mathtt{Choice}_\alpha^{\overline{w}}$ such that $\mathcal{M} ,v\nVdash \phi \mbox{ for some } v\in N$. By induction hypothesis, we get that $\mathcal{M}^T_W,\langle\overline{w},h_v\rangle\nvDash\phi$. We have that $N^T\in \mathbf{Choice}_\alpha^{\overline{w}}$. Therefore, our assumption entails that there exists $L_1^T\in \mathbf{Choice}_\alpha^{\overline{w}} \mbox{ such that } N^T\prec L_1^T \mbox{ and if } U=L_1^T\mbox{ or } L_1^T\preceq U, \mbox{ then } \mathcal{M} ^T,\langle\overline{w},h_u\rangle\vDash \phi \mbox{ for every } h_u\in U$. Here, we claim that $L_1$ is the choice cell that will witness to the fact that $\mathcal{M} ,w\Vdash \odot[\alpha ] \phi$. In order to show this, we first notice that the first item of \textbf{Claim 1} renders that $N\prec L_1$. Next, for every $v\in L_1$, $h_v\in L_1^T$ and therefore $\mathcal{M} ^T,\langle\overline{w},h_v\rangle\vDash \phi$; in this way, we have that for every $v\in L_1$, the induction hypothesis gives that $\mathcal{M} ,v\Vdash \phi$. Finally, let $L\in \mathtt{Choice}^{\overline{w}}_\alpha$ such that $L_1\preceq L$. By the first item of \textbf{Claim 1}, this last condition implies that $L_1^T\preceq L^T$. Take $v\in L$. We know that $h_v\in L^T$, so that the fact that $L_1^T\preceq L^T$ implies with our assumption that $\mathcal{M} ^T,\langle\overline{w},h_v\rangle\vDash \phi$; therefore, the induction hypothesis gives that $\mathcal{M} ,v\Vdash \phi$. With this, we have shown that our claim is true and that $\mathcal{M}_{w},w\Vdash \odot[\alpha ] \phi$.

\item (``$\odot_\mathcal{S}[\alpha ]$'')

($\Rightarrow$) We assume that $\mathcal{M} ,w\Vdash \odot_\mathcal{S}[\alpha ]\phi$. Then for every $L \in \mathtt{Choice}_\alpha^{\overline{w}}$ such that $\mathcal{M} ,v_L\nVdash \phi$ for some $v_L\in \sembrack{L}_\alpha^{\overline{w'}}$ (with $w\approx_\alpha w'$),
there exists $L'\in \mathtt{Choice}_\alpha^{\overline{w}}$ such that $ L\prec_s L'$ and if $ L''=L'$ or $ L'\preceq_s L''$, then for every $\overline{w''}$ such that $w\approx_\alpha w''$,  $\mathcal{M} ,v\Vdash \phi$ for every $v\in \sembrack{L''}_\alpha^{\overline{w''}}$. Let $N\in \mathbf{Choice}_\alpha^{\overline{w}}$ such that $\mathcal{M}^T_W,\langle\overline{w'},h_z\rangle\nvDash \phi \mbox{ for some } h_z\in [N]_\alpha^{\overline{w'}}$ (with $\overline{w}\sim_\alpha\overline{w'}$). By induction hypothesis, we get that $\mathcal{M} ,z\nVdash\phi$. The cell $N$ is $L_0^T$ for some $L_0\in\mathtt{Choice}^{\overline{w}}_\alpha$. Moreover, our definition of $\sim_\alpha$ in terms of $\approx_\alpha$ entails that $z\in \sembrack{L_0}_\alpha^{\overline{w'}}$, so that the fact that $\mathcal{M} ,z\nVdash\phi$ entails with our assumption that there exists $L_1\in \mathtt{Choice}_\alpha^{\overline{w}} \mbox{ such that } L_0\prec_s L_1 \mbox{ and if } L''=L_1\mbox{ or } L_1\preceq_s L''$ then for every $w''$ such that $w\approx_\alpha w''$, $\mathcal{M} ,v\Vdash \phi \mbox{ for every } v\in \sembrack{L''}_\alpha^{\overline{w''}}$. We claim that $L_1^T$ is the choice cell at moment $\overline{w}$ that will witness to the fact that $\mathcal{M} ^T,\langle\overline{w},h_w\rangle\vDash \odot_\mathcal{S}[\alpha ] \phi$. In order to show this, we first notice that the second item of \textbf{Claim 1}  renders that $L_0^T\prec_s L_1^T$. We also have that for every $\overline{w''}$ such that $\overline{w}\sim_\alpha\overline{w''}$ and $h_v\in [L_1^T]_\alpha^{\overline{w''}}$, $v\in \sembrack{L_1}_\alpha^{\overline{w''}}$ and therefore $\mathcal{M} ,v\Vdash \phi$; in this way, we have that for every $\overline{w''}$ such that $\overline{w}\sim_\alpha\overline{w''}$ and every $h_v\in [L_1^T]_\alpha^{\overline{w''}}$, the induction hypothesis gives that $\mathcal{M} ^T,(\overline{w''},h_v)\vDash \phi$. Finally, let $V^T\in \mathbf{Choice}^{\overline{w}}_\alpha$ such that $L_1^T\preceq_s V^T$. By the second item of \textbf{Claim 1} , this last condition implies that $L_1\preceq_s V$. Take $\overline{w''}$ such that $\overline{w}\sim_\alpha\overline{w''}$, and take $h_v\in [V^T]_\alpha^{\overline{w''}}$. We know that $v\in \sembrack{V}_\alpha^{\overline{w''}}$, so that the fact that $L_1\preceq_s V$ implies with our assumption that $\mathcal{M} ,v\Vdash \phi$; therefore, the induction hypothesis gives that $\mathcal{M} ^T,(\overline{w''},h_v)\vDash \phi$. With this, we have shown that our claim is true and that $\mathcal{M} ^T,\langle\overline{w},h_w\rangle\vDash \odot_{\mathcal{S}} [\alpha ] \phi$.       

($\Leftarrow$) We assume that $\mathcal{M} ^T,\langle\overline{w},h_w\rangle\vDash \odot_{\mathcal{S}}[\alpha ] \phi$. Then for every $L^T \in \mathbf{Choice}_\alpha^{\overline{w}}$ such that $\mathcal{M} ^T,\langle\overline{w},h_z\rangle\nvDash \phi$ for some $h_z\in [L^T]_\alpha^{\overline{w'}}$ (with $\overline{w}\sim_\alpha \overline{w'}$),
there exists $L^T_1\in \mathbf{Choice}_\alpha^{\overline{w}}$ such that $ L^T\prec_s L^T_1$ and if $U=L_1^T$ or $L^T\preceq_s U$, then for every $\overline{w''}$ such that $\overline{w}\sim_\alpha\overline{w''}$, $\mathcal{M} ^T,(\overline{w''},h_u)\vDash \phi$ for every $h_u\in [U]_\alpha^{\overline{w''}}$. Let $N\in \mathtt{Choice}_\alpha^{\overline{w}}$ such that $\mathcal{M} ,v\nVdash \phi \mbox{ for some } v\in \sembrack{N}_\alpha^{\overline{w'}}$ (with $w\approx_\alpha w'$). By induction hypothesis, we get that $\mathcal{M}^T,\langle\overline{w'},h_v\rangle\nvDash\phi$. The fact that $w\approx_\alpha w'$ implies that $\overline{w}\sim_\alpha\overline{w'}$. Since $h_v\in [N^T]_\alpha^{\overline{w'}}$, then, we know that there exists $L_1^T\in \mathbf{Choice}_\alpha^{\overline{w}}$ such that $ N^T\prec_s L_1^T$ and if $U=L_1^T$ or $L_1^T\preceq_s U$ ,then for every $\overline{w''}$ such that $\overline{w}\sim_\alpha\overline{w''},  \mathcal{M} ^T,(\overline{w''},h_u)\vDash \phi$ for every  $h_u\in [U]_\alpha^{\overline{w''}}$. Here, we claim that $L_1$ is the choice cell that will witness to the fact that $\mathcal{M} ,w\Vdash \odot_{\mathcal{S}}[\alpha ] \phi$. In order to show this, we first notice that the second item of \textbf{Claim 1} renders that $N\prec_s L_1$. Next, for every $w''$ such that $w\approx_\alpha w''$ (which implies that $\overline{w}\sim_\alpha \overline{w''})$ and every $v\in \sembrack{L_1}_\alpha^{\overline{w''}}$, $h_v\in [L_1^T]_\alpha^{\overline{w''}}$ and therefore $\mathcal{M} ^T,(\overline{w''},h_v)\vDash \phi$; in this way, we have that for every for every $w''$ such that $w\approx_\alpha w''$ and every $v\in \sembrack{L_1}_\alpha^{\overline{w''}}$, the induction hypothesis gives that $\mathcal{M} ,v\Vdash \phi$. Finally, let $L\in \mathtt{Choice}^{\overline{w}}_\alpha$ such that $L_1\preceq_s L$. By the second item of \textbf{Claim 1} , this last condition implies that $L_1^T\preceq_s L^T$. Take $w'$ such that $w\approx_\alpha w'$ (which implies that $\overline{w}\sim_\alpha \overline{w'})$, and take  $v\in \sembrack{L}_\alpha^{\overline{w}}$. We know that $h_v\in [L^T]_\alpha^{\overline{w''}}$, so that the fact that $L_1^T\preceq_s L^T$ implies with our assumption that $\mathcal{M} ^T,(\overline{w''},h_v)\vDash \phi$; therefore, the induction hypothesis gives that $\mathcal{M} ,v\Vdash \phi$. With this, we have shown that our claim is true and that $\mathcal{M} ,w\Vdash \odot_{\mathcal{S}}[\alpha ] \phi$.
\end{itemize}

\end{proof}

\subsection{Canonical Kripke-estit models}

We will show that the proof system $\Lambda$ is complete with respect to the class of  Kripke-estit models (and thus use Proposition \ref{avion1} to show completeness with respect to the class of epistemic utilitarian bi-valued BT-models). The strategy is to build a canonical structure from the syntax.  
\begin{definition}[Canonical Structure]
\label{dicaprio}
The tuple \[\mathcal{M}=\langle W^\Lambda, R_\square, \mathtt{Choice}, \{\mathtt{\approx}_\alpha\}_{\alpha\in Ags}, \mathtt{Value}_\mathcal{O}, \mathtt{Value}_\mathcal{S},\mathcal{V} 
\rangle \] is called a canonical structure for $\Lambda$ iff


\begin{itemize}
\item $W^\Lambda=\{w ;w  \mbox{ is a } \Lambda\textnormal{-MCS}\}$. 
\item $R_\square$ is a relation over $W^\Lambda$ defined by the following rule:  for $w,v\in W^\Lambda$, $wR_{\square}v$ iff for every $\phi$, $\square\phi\in w\Rightarrow \phi\in v$. For $w\in W^\Lambda$, the set $\{v\in W^\Lambda ;wR_\square v \}$ is denoted by $\overline{w}$. 
\item $\mathtt{Choice}$ is a function that assigns to each $\alpha$ and $\overline{w}$ a subset of $\mathcal{P}(\overline{w})$, which will be denoted by $\mathtt{Choice}_\alpha^{\overline{w}}$ and defined as follows: let $R_\alpha^{\overline{w}}$ be a relation on ${\overline{w}}$ such that for $w,v\in W^\Lambda$, $wR_{\alpha}^{\overline{w}}v$ iff for every $\phi$, $[\alpha ]\phi\in w\Rightarrow \phi\in v$; if we take $\mathtt{Choice}^{\overline{w}}_\alpha(v)=\{u\in \overline{w} ; vR_\alpha^{\overline{w}} u\}$, then we set $\mathtt{Choice}^{\overline{w}}_\alpha=\bigcup_{v\in \overline{w}}\mathtt{Choice}^{\overline{w}}_\alpha(v)$.

\item For each $\alpha \in Ags$, $\approx_\alpha$ is an epistemic relation on $W^\Lambda$ given by the following rule:  for $w,v\in W^\Lambda$, $w\approx_{\alpha}v$ iff for every $\phi$, $K_\alpha\phi\in w\Rightarrow \phi\in v$. 

\item $\mathtt{Value}_\mathcal{O}$ and $\mathtt{Value}_\mathcal{S}$ are defined as follows. For $\alpha\in Ags$ and $w\in W^\Lambda$, we first define two sets: \[\Sigma_\alpha^{w}=\{[\alpha ]\phi ;\odot [\alpha ]\phi \in w\}\]\[\Gamma_\alpha^{w}=\{K_\alpha\phi ;\odot_\mathcal{S}[\alpha ]\phi \in w\}.\] If we take $\Sigma^{w}=\bigcup_{\alpha\in Ags} \Sigma_\alpha^{w}$ and $\Gamma^{w}= \bigcup_{\alpha\in Ags} \Gamma_\alpha^{w}$, then we define the deontic functions by

\[\mathtt{Value} _\mathcal{O}(w) \ \ = \ \ \left\{\begin{array}{lll} 1 \mbox{ iff }  \Sigma^{w}  \subseteq w& \\ 0 \mbox{ otherwise}.&\end{array}\right.\]
\[\mathtt{Value} _\mathcal{S}(w) \ \ = \ \ \left\{\begin{array}{lll} 1 \mbox{ iff }  \Gamma^{w} \subseteq w& \\ 0 \mbox{ otherwise}.&\end{array}\right.\]

\item $\mathcal{V}$ is the canonical valuation, defined such that $w\in \mathcal{V}(p)$ iff $p\in w$. 

\end{itemize}
\end{definition}

\begin{proposition}\label{can}
The canonical structure $\mathcal{M}$ for $\Lambda$ is a Kripke-estit model.
\end{proposition}

\begin{proof}
We need to show that the tuple $\langle W^\Lambda, R_\square, \mathtt{Choice}, \{\mathtt{\approx}_\alpha\}_{\alpha\in Ags}, \mathtt{Value}_\mathcal{O}, \mathtt{Value}_\mathcal{S}\rangle$ is a Kripke-estit frame, which amounts to showing that the tuple validates the five items in the definition of Krioke-estir models. 

\begin{enumerate}

\item It is clear that $R_\square$ is an equivalence relation, since $\Lambda$ includes the $\mathbf{S5}$ axioms for $\square$. 
\item Since $\Lambda$ includes the $\mathbf{S5}$ axioms for $[\alpha ]$ (for each $\alpha\in Ags$), we have that $R_\alpha^{\overline{w}}$ is an equivalence relation for each $\alpha\in Ags$ and $w\in W^\Lambda$. Moreover, since $\Lambda$ includes $\square p\to [\alpha ]p$ as an axiom schema, we have that $R_\alpha\subseteq R_\square $. In this way, $\mathtt{Choice}$ indeed assigns to each $\alpha$ and $\overline{w}$ a partition of $\overline{w}$. Now, we must verify that $\mathcal{M}$ validates the constraint $(\mathtt{IA})_K$. In order to do that, we need two intermediate results: 
\begin{enumerate}[a)]
\item For a fixed $w_*\in W^\Lambda$, we have that $w\in\overline{w_*}$ iff $\{\square \psi; \square\psi\in w_*\}\subseteq w$. 

$(\Rightarrow)$ Let $w\in\overline{w_*}$ (which means that $w_*R_\square w$). Take $\phi$ a formula of $\mathcal{L}_{\textsf{KO}}$ such that $\square\phi\in w_*$. Since $w_*$ is closed under \emph{Modus Ponens}, the $(4)$ axiom for $\square$ implies that $\square\square\phi\in w_*$ as well. Therefore, by the definition of $R_\square$, we get that $\square\phi \in w$. 

$(\Leftarrow)$ We assume that $\{\square \psi; \square\psi\in w_*\}\subseteq w$. Take $\phi$ a formula of $\mathcal{L}_{\textsf{KO}}$ such that $\square\phi\in w_*$. By our assumption, we get that $\square \phi\in W$. Since $w$ is closed under \emph{Modus Ponens}, the $(T)$ axiom for $\square$ implies that $\phi\in w$ as well. In this way, we have that the fact that $\square \phi\in w_*$ implies that $\phi\in w$, which means that $w_*R_\square w$ and $w\in\overline{w_*}$. 

\item For a fixed $w_*\in W^\Lambda$ and $s:Ags\to \mathcal{P}(\overline{w_*})$ a function that maps $\alpha$ to a member of $\mathtt{Choice}^{\overline{w_*}}_\alpha$ such that $v_{s(\alpha)}\in s(\alpha)$, we have that $w\in s(\alpha)$ iff $\Delta_{s(\alpha)}=\{[\alpha ] \psi ;[\alpha ]\psi\in v_{s(\alpha)}\}\subseteq w$.

$(\Rightarrow)$ Let $w\in s(\alpha)$ (which means that $v_{s(\alpha)}R_\alpha^{\overline{w_*}} w$). Take $\phi$ a formula of $\mathcal{L}_{\textsf{KO}}$ such that $[\alpha ]\phi\in v_{s(\alpha)} $. Since $v_{s(\alpha)}$ is closed under \emph{Modus Ponens}, the $(4)$ axiom for $[\alpha ]$ implies that $[\alpha ][\alpha ]\phi\in v_{s(\alpha)}$ as well. Therefore, by definition of $R_\alpha^{\overline{w_*}}$, we get that $[\alpha ]\phi \in w$. 

$(\Leftarrow)$ We assume that $\Delta_{s(\alpha)}=\{[\alpha ] \psi ; [\alpha ]\psi\in v_{s(\alpha)}\}\subseteq w$. Take $\phi$ a formula of $\mathcal{L}_{\textsf{KO}}$ such that $[\alpha ]\phi\in v_{s(\alpha)}$. By our assumption, we get that $[\alpha ]\phi\in w$. Since $w$ is closed under \emph{Modus Ponens}, the $(T)$ axiom for $[\alpha ]$ implies that $\phi\in w$ as well. In this way, we have that the fact that $[\alpha ] \phi\in v_{s(\alpha)}$ implies that $\phi\in w$, which means that $v_{s(\alpha)}R_\alpha^{\overline{w_*}} w$ and $w\in s(\alpha)$.
\end{enumerate}

Next, we will show that for a fixed $w_*\in W^\Lambda$ and $s:Ags\to \mathcal{P}(\overline{w_*})$ a function that maps $\alpha$ to a member of $\mathtt{Choice}^{\overline{w_*}}_\alpha$ such as in item b) above, we have that $\bigcup_{\alpha \in Ags}\Delta_{s(\alpha)}\cup \{\square \psi ; \square\psi\in w_*\}$ is $\Lambda$-consistent.

First, we will show that $\bigcup_{\alpha\in Ags} \Delta_{s(\alpha)}$ is consistent. Suppose that this is not the case. Then there exists a set $\{\phi_1,\dots,\phi_n\}$ of formulas of $\mathcal{L}_{\textsf{KO}}$ such that $[\alpha_{i} ]\phi_i\in v_{s(\alpha_i)}$ for every $1\leq i \leq n$ and \begin{equation}\label{paul}
\vdash([\alpha_{1} ]\phi_1\wedge\dots\wedge [\alpha_{n}  ]\phi_n)\to \bot.
\end{equation}

Without loss of generality, we assume that $\alpha_i\neq\alpha_j$ for all $j\neq i$ such that $j,i\in \{1,\dots,n\}$ --this assumption hinges on the fact that any stit operator distributes over conjunction. Notice that the fact that $[\alpha_{i} ]\phi_i\in v_{s(\alpha_i)}$ for every $1\leq i \leq n$ implies that $\Diamond[\alpha_{i} ] \in w_*$ for every $1\leq i \leq n$. 
Since $w_*$ is closed under conjunction, we also have that $\Diamond [\alpha_{1} ]\phi_1\wedge\dots\wedge\Diamond[\alpha_{n} ]\phi_n\in w_*$.

By the independence-of-agents axiom, we have that \begin{equation}\label{paul2} 
\vdash\Diamond [\alpha_{1} ]\phi_1\wedge\dots\wedge\Diamond[\alpha_{n} ]\phi_n\to \Diamond( [\alpha_{1} ]\phi_1\wedge\dots\wedge[\alpha_{n} ]\phi_n).
\end{equation} Therefore, equations \eqref{paul2} and \eqref{paul},  imply that \begin{equation}\label{paul3} 
\vdash\Diamond [\alpha_{1} ]\phi_1\wedge\dots\wedge\Diamond[\alpha_{n} ]\phi_n\to \Diamond \bot.
\end{equation}
But this is a contradiction, since we had seen that $\Diamond [\alpha_{1} ]\phi_1\wedge\dots\wedge\Diamond[\alpha_{n} ]\phi_n\in w_*$, and $w_*$ is a $\Lambda$-MCS. Therefore, $\bigcup_{\alpha\in Ags} \Delta_{s(\alpha)}$ is consistent.

Next, we show that the union $\bigcup_{\alpha\in Ags} \Delta_{s(\alpha)}\cup  \{\square\psi ; \square\psi \in w_*\}$ is also consistent. Suppose that this is not the case. Since $\bigcup_{\alpha\in Ags} \Delta_{s(\alpha)}$ is consistent, there must exist sets $\{\phi_1,\dots,\phi_n\}$ and $\{\theta_1,\dots,\theta_m\}$ of formulas of $\mathcal{L}_{\textsf{KO}}$ such that $[\alpha_i ]\phi_i\in v_{s(\alpha_i)}$ for every $1\leq i \leq n$, $\square\theta_i\in w_*$ for every $1\leq i\leq m$, and \begin{equation}\label{newman} \vdash([\alpha_1 ]\phi_1\wedge\dots\wedge [\alpha_n]\phi_n)\wedge (\square\theta_1\wedge\dots\wedge\square\theta_m)\to\bot. 
\end{equation} Let $\theta=\theta_1\land\dots \land\theta_m$. Since $\square$ distributes over conjunction, we have that  $\square\theta=\square\theta_1\wedge\dots\wedge\square\theta_m$, where it is important to mention that since $w_*$ is a $\Lambda$-MCS, then $\square\theta\in w_*$. In these terms, \eqref{newman} implies that \begin{equation}\label{newman2}\vdash([\alpha_1 ]\phi_1\wedge\dots\wedge [\alpha_n ]\phi_n)\to\lnot \square \theta.\end{equation} Once again, we assume without loss of generality that $\alpha_i\neq\alpha_j$ for all $j\neq i$ such that $j,i\in \{1,\dots,n\}$. Analogous to the procedure we used to show that $\bigcup_{\alpha\in Ags} \Delta_{s(\alpha)}$ is consistent, \eqref{newman2} implies that \begin{equation}\label{newman3} 
\vdash\Diamond [\alpha_{1} ]\phi_1\wedge\dots\wedge\Diamond[\alpha_{n} ]\phi_n\to \Diamond \lnot \square\theta.
\end{equation}

This entails that $\Diamond \lnot \square\theta\in w_*$, but this is a contradiction, since the fact that $\square\theta\in w_*$ implies with the $(4)$ axiom for $\square$ that $\square\square\theta\in w_*$. 

Now, let $u_*$ be the $\Lambda$-MCS that includes $\bigcup_{\alpha\in Ags} \Delta_{s(\alpha)}\cup  \{\square\psi ; \square\psi \in w_*\}$. By the intermediate result  a), it is clear that $u_*\in \overline{w_*}$. By the intermediate result b), it is clear that $u_*\in s(\alpha)$ for every $\alpha\in Ags$. Therefore, we have shown that for a fixed $w_*\in W$, we have that each function $s:Ags\to \mathcal{P}(\overline{w_*})$ that maps $\alpha$ to a member of $\mathtt{Choice}^{\overline{w_*}}_\alpha$ is such that $\bigcap_{\alpha \in Ags} s(\alpha) \neq \emptyset$, which means that $\mathcal{M}$ validates the constraint $(\mathtt{IA})_K$.

\item Since the axiom system $\Lambda$ includes the $\mathbf{S5}$ axioms for $K_\alpha$ for each $\alpha\in Ags$, we have that $\approx_\alpha$ is an equivalence relation for each $\alpha\in Ags$. We must now verify that $\mathcal{M}$ validates the constraints $(\mathtt{OAC})_K$ and $(\mathtt{Unif-H})_K$.
\begin{enumerate}[i)]
\item For $(\mathtt{OAC})_K$, fix $w_*\in W^\Lambda$ and $\alpha\in Ags$. We assume that $v,u\in \overline{w_*}$ are such that $v\approx_\alpha u$. Let $v'\in \mathtt{Choice}^{\overline{w_*}}_\alpha(v)$. This means that $vR_\alpha v'$. We want to show that $v'\approx_\alpha u$, so let $\phi$ be a formula of $\mathcal{L}_{\textsf{KO}}$ such that $K_\alpha \phi\in v'$. By the $(4)$ axiom for $K_\alpha$, we have that $K_\alpha K_\alpha\phi \in v'$. Similarly, since all substitutions of the $(OAC)$ axiom lie within $v'$ and it is closed under \emph{Modus Ponens}, we get that $[\alpha ]K_\alpha\phi$ also lies in $v'$. Since $v'R_\alpha v$, this implies that $K_\alpha\phi\in v$. Therefore, our assumption that $v\approx_\alpha u$ entails that $\phi \in u$. With this, we have shown that the fact that $K_\alpha\phi\in v'$ implies that $\phi\in u$, which means that $v'\approx_\alpha u$.

\item For $(\mathtt{Unif-H})_K$, fix $w_1,w_2\in W^\Lambda$. We assume that there exist $v\in \overline{w_1}$ and $u\in\overline{w_2}$ such that $v\approx_\alpha u$. Take $v'\in\overline{w_1}$.

We will show that $u''=\{\psi ; K_\alpha\psi\in v'\}\cup \{\square\psi ; \square\psi\in u \}$ is consistent. In order to do so, we will first show that $\{\psi ; K_\alpha\psi\in v'\}$ is consistent. Suppose for a contradiction that it is not consistent. Then there exists a set $\{\psi_1,\dots,\psi_n\}$ of formulas of $\mathcal{L}_{\textsf{KO}}$ such that $\{\psi_1,\dots,\psi_n\}\subseteq \{\psi ; K_\alpha\psi\in v'\}$ and $\vdash \psi_1\wedge\dots\wedge \psi_n\to \bot$ (a). Now, the fact that $\{\psi_1,\dots,\psi_n\}\subseteq \{\psi ; K_\alpha\psi\in v'\}$ means that $K_\alpha \psi_i\in v'$ for every $1\leq i \leq n$; by necessitation of $K_\alpha$ and its distributivity over conjunction, we get that (a) implies that $\vdash K_\alpha\psi_1\wedge\dots\wedge K_\alpha\psi_n\to K_\alpha\bot$, but this is a contradiction, since $
v'$ is a $\Lambda$-MCS which, being, includes $K_\alpha\psi_1\wedge\dots\wedge K_\alpha\psi_n$. 

Next, we show that $u''=\{\psi ; K_\alpha\psi\in v'\}\cup \{\square\psi ; \square\psi\in u \}$ is also consistent. Suppose for a contradiction that it is not consistent. Since $\{\psi ; K_\alpha\psi\in v'\}$ is consistent, there must exist sets $\{\phi_1,\dots,\phi_n\}$ and $\{\theta_1,\dots,\theta_m\}$ of formulas of $\mathcal{L}_{\textsf{KO}}$ such that $\{\psi_1,\dots,\psi_n\}\subseteq \{\psi ; K_\alpha\psi\in v'\}$, $\square\theta_i\in w_2$ for every $1\leq i\leq m$, and $\vdash(\phi_1\wedge\dots\wedge\phi_n)\wedge (\square\theta_1\wedge\dots\wedge\square\theta_m)\to\bot$ (a). Let $\theta=\theta_1\land\dots \land\theta_m$  and $\phi=\phi_1\land\dots \land\phi_n$. Since $\square$ distributes over conjunction, we have that $\square\theta=\square\theta_1\wedge\dots\wedge\square\theta_m$, where it is important to mention that since $u$ is a $\Lambda$-MCS, then $\square\theta\in u$ and $\square \square\theta\in u$ as well ($\star$). In this way, (a) implies that $\vdash\phi\to\lnot \square \theta$ and thus that $\vdash\Diamond\phi\to\Diamond\lnot\square\theta$ (b). Now, by necessitation of $K_\alpha$, we get that (b) implies that $\vdash K_\alpha\Diamond \phi \to K_\alpha\Diamond\lnot \square \theta$ (c). Notice that the facts that $K_\alpha\phi_i\in v'$ for every $1\leq i\leq n$, that $K_\alpha$ distributes over conjunction, and that $v'$ is a $\Lambda$-MCS imply that $K_\alpha\phi=K_\alpha\phi_1\wedge\dots\wedge K_\alpha\phi_n\in v'$. The fact that $v'\in \overline {w_1}=\overline{v}$ implies that $\Diamond K_\alpha \phi \in v$, so that (c) and the facts that $v$ contains all substitutions of the $(Unif-H)$ axiom and is closed under \emph{Modus Ponens} entail that $K_\alpha\Diamond\phi\in v$. Now, this last inclusion implies, with our assumption that $v\approx_\alpha u$, that $\Diamond\phi\in u$, which by (b) in turn yields that $\Diamond\lnot\square\theta\in u$, contradicting $(\star)$. Therefore, $u''$ is consistent. 

Finally, let $u'$ be the $\Lambda$-MCS that includes $u''$. It is clear from its construction that $u'\in\overline{u}=\overline{w_2}$ and that $v'\approx_\alpha u'$, 

With this, we have shown that $\mathcal{M}$ validates the constraint $(\mathtt{Unif-H})_K$.

\end{enumerate}
\item Both $\mathtt{Value}_{\mathcal{O}}$ and $\mathtt{Value}_{\mathcal{S}}$ are well defined functions with range within $\mathds{R}$.
\end{enumerate}

\end{proof}

As is usual with canonical structures, our objective is to prove the so-called \emph{truth} lemma, which says that for every formula $\phi$ of $\mathcal{L}_{\textsf{KO}}$ and every $w\in W^\Lambda$, we have that $
\mathcal{M}^n ,w\Vdash \phi \ \textnormal{iff} \ \phi\in w.$ This is done by induction on $\phi$, and the inductive step for each modal operator requires previous results (such as the important \emph{existence} lemmas). In the case of $\square, K_\alpha$, and $[\alpha]$, these previous results are standard (Lemma \ref{carajo} below). However, in the case of $\odot[\alpha]$ and $\odot_{\mathcal{S}}[\alpha]$, the previous results require a lot of work. In Lemmas \ref{vang} and \ref{kid} below, we pinpoint the main claims that we need to show the inductive step of the truth lemma for the ought operators. Each item in these two lemmas is shown by a series of sub-lemmas, and the reader is advised to carefully go over the full proofs provided in the supplement to the present work.     

\begin{lemma}[Existence] \label{carajo} 
Let $\mathcal{M}$ be the canonical Kripke-estit model for $\Lambda$. Let $w\in W^\Lambda$. For a given formula $\phi$ of $\mathcal L_{\textsf{KO}}$, the following hold: 
\begin{enumerate}
\item  $\square \phi \in w$ iff $\phi\in v$ for every $v\in \overline{w}$. 
\item $[\alpha ]\phi\in w$ iff $\phi\in v$ for every $v\in \overline{w}$ such that $wR_\alpha v$.
\item $K_\alpha\phi\in w$ iff $\phi\in v$ for every $v\in W^\Lambda$ such that $w\approx_\alpha v$.
\end{enumerate}
\end{lemma} 
\begin{proof}
Let $w\in W^\Lambda$, and take $\phi$ a formula of $\mathcal L_{\textsf{KO}}$.

\begin{enumerate}

\item $(\Rightarrow)$ We assume that $\square \phi \in w$. If $v$ lies within $\overline{w}$, this means that $wR_\square v$ holds. Therefore, $\phi\in v$. 

$(\Leftarrow)$ We work by contraposition. Assume that $\square\phi\notin w$. We will show that there is a world in $\overline{w}$ such that $\phi$ does not lie within it. For this, let $v'=\{\psi ; \square\psi\in w\}$, which is consistent by virtue of the following argument: suppose for a contradiction that $v'$ is not consistent; then there exists a set $\{\psi_1,\dots,\psi_n\}$ of formulas of $\mathcal{L}_{\textsf{KO}}$ such that $\{\psi_1,\dots,\psi_n\}\subseteq v'$ and $\vdash \psi_1\wedge\dots\wedge \psi_n\to \bot$ (a); now, the fact that $\{\psi_1,\dots,\psi_n\}\subseteq v'$ means that $\square \psi_i\in w$ for every $1\leq i \leq n$; due to necessitation of $\square$ and its distributivity over conjunction, we get that (a) implies that $\vdash \square\psi_1\wedge\dots\wedge \square\psi_n\to \square\bot$, but this is a contradiction, since $w$ is a $\Lambda$-MCS which, being closed under conjunction, includes $\square\psi_1\wedge\dots\wedge \square\psi_n$. Now, we define $v''=v'\cup \{\lnot\phi\}$, and show that it is also consistent as follows: suppose for a contradiction that it is not consistent; since $v'$ is consistent, we have that there exists a set $\{\psi_1,\dots,\psi_n\}$ of formulas of $\mathcal{L}_{\textsf{KO}}$ such that $\{\psi_1,\dots,\psi_n\}\subseteq v'$ and $\vdash \psi_1\wedge\dots\wedge \psi_n \wedge \lnot\phi\to \bot$, which then implies that $\vdash \psi_1\wedge\dots\wedge \psi_n \to \phi$ (b); due to necessitation of $\square$ and its distributivity over conjunction, we get that (b) implies that $\vdash \square\psi_1\wedge\dots\wedge \square\psi_n\to \square\phi$ (c); but notice that since $w$ is a $\Lambda$-MCS, then $\square\psi_1\wedge\dots\wedge \square\psi_n\in w$, so that (c) and the fact that $w$ is closed under \emph{Modus Ponens} entail that $\square\phi\in w$, contradicting the initial assumption that $\square\phi\notin w$. Finally, let $v$ be the $\Lambda$-MCS that includes $v''$. On one hand, it is clear from its construction that $\phi \notin v$ and that $wR_\square v$. Thus, we have shown that assuming that $\square\phi\notin w$ implies the existence of a $v\in \overline{w}$  such that $\phi\notin v$.


\item $(\Rightarrow)$ We assume that $[\alpha ]\phi\in w$. Let $v\in \overline{w}$ such that $wR_\alpha v$. The definition of $R_\alpha$ straightforwardly gives that $\phi\in v$.  

$(\Leftarrow)$ We work by contraposition. Assume that $[\alpha ]\phi\notin w$. We will show that there is a world $v$ in $\overline{w}$ such that $wR_\alpha v$ and such that $\phi$ does not lie within it. For this, let $v'=\{\psi ; [\alpha ]\psi\in w\}$, which is consistent by a similar argument than the one introduced in the proof of the above item: suppose for a contradiction that $v'$ is not consistent; then there exists a set $\{\psi_1,\dots,\psi_n\}$ of formulas of $\mathcal{L}_{\textsf{KO}}$ such that $\{\psi_1,\dots,\psi_n\}\subseteq v'$ and $\vdash \psi_1\wedge\dots\wedge \psi_n\to \bot$ (a); now, the fact that $\{\psi_1,\dots,\psi_n\}\subseteq v'$ means that $[\alpha ] \psi_i\in w$ for every $1\leq i \leq n$; necessitation of $[\alpha ]$ and its distributivity over conjunction yield that (a) implies that $\vdash [\alpha ]\psi_1\wedge\dots\wedge [\alpha ]\psi_n\to [\alpha ]\bot$, but this is a contradiction, since $w$ is a $\Lambda$-MCS which, being closed under conjunction, includes $[\alpha ]\psi_1\wedge\dots\wedge [\alpha ]\psi_n$. Now, we define $v''=v'\cup \{\lnot\phi\}$, and we show that it is also consistent as follows: suppose for a contradiction that it is not consistent; since $v'$ is consistent, we have that there exists a set $\{\psi_1,\dots,\psi_n\}$ of formulas of $\mathcal{L}_{\textsf{KO}}$ such that $\{\psi_1,\dots,\psi_n\}\subseteq v'$ and $\vdash \psi_1\wedge\dots\wedge \psi_n \wedge \lnot\phi\to \bot$, which then implies that $\vdash \psi_1\wedge\dots\wedge \psi_n \to \phi$ (b); due to necessitation of $[\alpha ]$ and its distributivity over conjunction, we get that (b) implies that $\vdash [\alpha ]\psi_1\wedge\dots\wedge [\alpha ]\psi_n\to [\alpha ]\phi$ (b); but notice that since $w$ is a $\Lambda$-MCS, then $[\alpha ]\psi_1\wedge\dots\wedge [\alpha ]\psi_n\in w$, so that (b) and the fact that $w$ is closed under \emph{Modus Ponens} entail that $[\alpha ]\phi\in w$, contradicting the initial assumption that $[\alpha ]\phi\notin w$. Finally, let $v$ be the $\Lambda$-MCS that includes $v''$. On one hand, it is clear from its construction that $\phi \notin v$ and that $wR_\alpha v$, by definition of $R_\alpha$. On the other, axiom $(A2)$ renders that actually $v\in \overline{w}$ (if $\square\theta \in w$, then $[\alpha ]\theta \in w$ and thus $ \theta \in v$). Therefore, we have shown that assuming that $[\alpha ]\phi\notin w$ implies the existence of a $v\in \overline{w}$ such that $wR_\alpha v$ and $\phi\notin v$.

\item $(\Rightarrow)$ We assume that $K_\alpha\phi\in w$. Let $v\in W^\Lambda$ such that $w\approx_\alpha v$. The definition of $\approx_\alpha$ straightforwardly gives that $\phi\in v$.  

$(\Leftarrow)$ We work by contraposition. Assume that $K_\alpha\phi\notin w$. We will show that there is a world $v$ in $W^\Lambda$ such that $w\approx_\alpha v$ and such that $\phi$ does not lie within it. For this, let $v'=\{\psi ; K_\alpha\psi\in w\}$, which is consistent by an argument similar to the one introduced in the proofs of the above items: suppose for a contradiction that $v'$ is not consistent; then there exists a set $\{\psi_1,\dots,\psi_n\}$ of formulas of $\mathcal{L}_{\textsf{KO}}$ such that $\{\psi_1,\dots,\psi_n\}\subseteq v'$ and $\vdash \psi_1\wedge\dots\wedge \psi_n\to \bot$ (a); now, the fact that $\{\psi_1,\dots,\psi_n\}\subseteq v'$ means that $K_\alpha \psi_i\in w$ for every $1\leq i \leq n$; necessitation of $K_\alpha$ and its distributivity over conjunction yield that (a) implies that $\vdash K_\alpha\psi_1\wedge\dots\wedge K_\alpha\psi_n\to K_\alpha\bot$, but this is a contradiction, since $w$ is a $\Lambda$-MCS which, being closed under conjunction, includes $K_\alpha\psi_1\wedge\dots\wedge K_\alpha\psi_n$. Now, we define $v''=v'\cup \{\lnot\phi\}$ $v''=v'\cup \{\lnot\phi\}$, and we show that it is also consistent as follows: suppose for a contradiction that it is not consistent; since $v'$ is consistent, we have that there exists a set $\{\psi_1,\dots,\psi_n\}$ of formulas of $\mathcal{L}_{\textsf{KO}}$ such that $\{\psi_1,\dots,\psi_n\}\subseteq v'$ and $\vdash \psi_1\wedge\dots\wedge \psi_n \wedge \lnot\phi\to \bot$, which then implies that $\vdash \psi_1\wedge\dots\wedge \psi_n \to \phi$ (b); due to necessitation of $K_\alpha$ and its distributivity over conjunction, we get that (b) implies that $\vdash K_\alpha \psi_1\wedge\dots\wedge K_\alpha \psi_n\to K_\alpha \phi$ (b); but notice that since $w$ is a $\Lambda$-MCS, then $K_\alpha \psi_1\wedge\dots\wedge K_\alpha \psi_n\in w$, so that (b) and the fact that $w$ is closed under \emph{Modus Ponens} entail that $K_\alpha\phi\in w$, contradicting the initial assumption that $K_\alpha \phi\notin w$. Finally, let $v$ be the $\Lambda$-MCS that includes $v''$. On one hand, it is clear from its construction that $\phi \notin v$ and that $w\approx_\alpha v$, by definition of $\approx_\alpha$. Thus, we have shown that assuming that $K_\alpha\phi\notin w$ implies the existence of a $v\in W^\Lambda$  such that $w\approx_\alpha v$ and $\phi\notin v$.
\end{enumerate} 
\end{proof}

\begin{lemma}[Existence for objective ought]\label{vang} Fix $\alpha\in Ags$ and $w\in W^\Lambda$. \begin{enumerate}[(a)]
\item For every formula $\phi$ of $\mathcal{L}_{\textsf{KO}}$, $\odot [\alpha ]\phi\in w$ iff $[\alpha ]\phi\in v$ for every $v\in \overline{w}\subseteq W^\Lambda$ such that $\Sigma_\alpha^{w}\subseteq v$.
\item For every $v\in \mathtt{Choice}^{\overline{w}}_\alpha(w)$, $\Sigma_\alpha^{w}\subseteq w$ iff  $\Sigma_\alpha^{w}\subseteq v$ 
\item $\Sigma_\alpha^{w}\subseteq w$ iff $\mathtt{Choice}^{\overline{w}}_\alpha(w)\in \mathtt{Optimal}^{\overline{w}}_\alpha$. 
\item For each $L\in \mathtt{Choice}^{\overline{w}}_\alpha-\mathtt{Optimal}^{\overline{w}}_\alpha$, there exists $L'\in \mathtt{Optimal}^{\overline{w}}_\alpha$ such that $L\prec L'$. Thus, $\mathcal{M} ,w\Vdash \odot[\alpha ]\phi  \textnormal{ iff }  \forall L \in \mathtt{Optimal}^{\overline{w}}_\alpha, \mathcal{M},v\Vdash \phi \mbox{ for every } v\in L.$
\end{enumerate}
 
\end{lemma}

\begin{proof}
\begin{enumerate}[(a)]
    \item \textbf{Claim 1.} Let $\mathcal{M}$ be the canonical Kripke-estit model for $\Lambda$. We have that the following hold:
\begin{itemize}
\item For $w_*\in W^\Lambda$, $w\in\overline{w_*}$ iff $\{\square \psi ; \square\psi\in w_*\}\subseteq w$. 
\item For any pair of agents $\alpha, \alpha\in Ags$ and any formula $\phi$ of $\mathcal{L}_{\textsf{KO}}$, $\odot [\alpha]\phi\to \lnot \odot [\alpha]\lnot\phi\in w$ for every $w\in W^\Lambda$. We will refer to this property as ``consistency of ought-to-do''. 
\end{itemize}
The proof of the first item goes as follows: $(\Rightarrow)$ Let $w\in\overline{w_*}$ (which means that $w_*R_\square w$). Take $\phi$ a formula of $\mathcal{L}_{\textsf{KO}}$ such that $\square\phi\in w_*$. Since $w_*$ is closed under \emph{Modus Ponens}, the $(4)$ axiom for $\square$ implies that $\square\square\phi\in w_*$ as well. Therefore, by definition of $R_\square$, we get that $\square\phi \in w$. $(\Leftarrow)$ We assume that $\{\square \psi ; \square\psi\in w_*\}\subseteq w$. Take $\phi$ a formula of $\mathcal{L}_{\textsf{KO}}$ such that $\square\phi\in w_*$. By our assumption, we get that $\square \phi\in W$. Since $w$ is closed under \emph{Modus Ponens}, the $(T)$ axiom for $\square$ implies that $\phi\in w$ as well. In this way, we have that the fact that $\square \phi\in w_*$ implies that $\phi\in w$, which means that $w_*R_\square w$ and $w\in\overline{w_*}$. 

The proof of the second item is a direct consequence from the fact that, for $\alpha,\beta\in Ags$ and $\phi$ a formula of $\mathcal{L}_{\textsf{KO}}$, we have that $\vdash\odot [\alpha ]\phi\to \lnot \odot [\beta ]\lnot\phi$. This comes from the fact that, for $\alpha,\beta\in Ags$, $\vdash\odot [\alpha ]\phi\land \odot [\beta]\lnot\phi\to \bot$ --which can be seen by applying axioms $(Oic)$, $(IA)$, and axiom $(T)$ for $[\alpha ]$ and $[\beta]$. 

Now we proceed with the proof of the main statement.

$(\Rightarrow)$ The definition of $\Sigma_\alpha^{w}$ implies straightforwardly that $[\alpha]\phi\in v$.

$(\Leftarrow)$ We work by contraposition. Suppose that $\odot [\alpha]\phi\notin w$. 

\begin{itemize}
\item We first show that $\Sigma_\alpha^{w}$ is consistent. 
Suppose that $\Sigma_\alpha^{w}$ is not consistent. Then there is a set $\{\phi_1,\dots,\phi_n\}$ of formulas of $\mathcal{L}_{\textsf{KO}}$ such that $\odot [\alpha]\phi_i\in w$ for every $1\leq i \leq n$ and  $\vdash([\alpha]\phi_1\wedge\dots\wedge [\alpha]\phi_n)\to \bot$. This last thing  implies that $\vdash([\alpha]\phi_1\wedge\dots\wedge [\alpha]\phi_{n-1})\to \lnot [\alpha]\phi_{n}$, so that if we take $\psi=\phi_1\wedge\dots\wedge \phi_{n-1}$, then distributivity of $[\alpha]$ over conjunction gives us that $\vdash[\alpha]\psi\to \lnot [\alpha]\phi_{n}$. By necessitation of the ought-to-do operator and its axiom $(K)$, we get that \begin{equation}\label{ay}\vdash\odot [\alpha]([\alpha]\psi)\to \odot [\alpha](\lnot [\alpha]\phi_{n}).
\end{equation} 
First, we notice that since $w$ is a $\Lambda$-MCS, then $\odot[\alpha]\psi\in w$ (by distributivity of $\odot[\alpha]$ over conjunction). The fact that $\odot[\alpha]\psi\in w$ implies with axiom $(A4)$ that $\odot [\alpha]([\alpha]\psi)\in w$, which with equation \eqref{ay} implies that $\odot [\alpha](\lnot [\alpha]\phi_{n})\in w$. However, the fact that $\odot[\alpha]\phi_n\in w$ implies that $\odot [\alpha]([\alpha]\phi_n)\in w$. Thus, we both have that $\odot [\alpha](\lnot [\alpha]\phi_{n})\in w$ and  that $\odot [\alpha]([\alpha]\phi_n)\in w$, which is a contradiction (by consistency of ought-to-do). 

\item Next, we observe that the union $\Sigma_\alpha^{w} \cup  \{\square\psi ; \square\psi \in w\}$ is also consistent. Suppose that $\Sigma_\alpha^{w} \cup  \{\square\psi ; \square\psi \in w\}$ is not consistent. Since $\Sigma_\alpha^{w}$ is consistent, we have that there must exist sets $\{\phi_1,\dots,\phi_n\}$ and $\{\theta_1,\dots,\theta_m\}$ of formulas of $\mathcal{L}_{\textsf{KO}}$ such that $\odot [\alpha]\phi_i\in w$ for every $1\leq i \leq n$, $\square\theta_i\in w$ for every $1\leq i\leq m$, and \begin{equation}\label{sh} \vdash([\alpha]\phi_1\wedge\dots\wedge [\alpha]\phi_n)\wedge (\square\theta_1\wedge\dots\wedge\square\theta_m)\to\bot. 
\end{equation} Let $\theta=\theta_1\land\dots \land\theta_m$. Since $\square$ distributes over conjunction, we have that  $\square\theta=\square\theta_1\wedge\dots\wedge\square\theta_m$, where it is important to mention that since $w$ is a $\Lambda$-MCS, then $\square\theta\in w$. In these terms, \eqref{sh} implies that $\vdash([\alpha]\phi_1\wedge\dots\wedge [\alpha]\phi_n)\to\lnot \square \theta $. By necessitation of the ought-to-do operator, its axiom $(K)$, and its distributivity over conjunction, we get that \begin{equation}\label{ay3}\vdash(\odot [\alpha]([\alpha]\phi_1)\land\dots \land \odot [\alpha]([\alpha]\phi_n)) \to \odot [\alpha](\lnot \square\theta).\end{equation} We have that $\odot [\alpha]([\alpha]\phi_i)\in w$ for every $1\leq i\leq n$. Therefore, their conjunction also lies in $w$, which by \eqref{ay3} entails that $\odot [\alpha](\lnot \square\theta)\in w$. However, we had also observed that $\square \theta \in w$, which by axiom $(4)$ for $\square$ implies that $\square \square \theta\in w$ and therefore that $\odot [\alpha] (\square \theta)\in w$. We have struck a contradiction, since we cannot have that both $\odot [\alpha](\lnot \square\theta)$ and $\odot [\alpha] (\square \theta)$ lie in $w$ (consistency of ought-to-do). Therefore, $\Sigma_\alpha^{w} \cup  \{\square\psi ; \square\psi \in w\}$ is consistent. 
\item Finally, we observe that the union $\Sigma_\alpha^{w} \cup \{\square\psi ; \square\psi \in w\}\cup \{\lnot [\alpha] \phi\}$ is also consistent:

First, recall that we are assuming that $\odot [\alpha]\phi\notin w$. Now, suppose that $\Sigma_\alpha^{w} \cup \{\square\psi ; \square\psi \in w\}\cup \{\lnot [\alpha] \phi\}$ is not consistent. Since $\Sigma_\alpha^{w}\cup \{\square\psi ; \square\psi \in w\}$ is consistent, then there must exist sets $\{\phi_1,\dots,\phi_n\}$ and $\{\theta_1,\dots,\theta_m\}$ of formulas of $\mathcal{L}_{\textsf{KO}}$ such that $\odot [\alpha]\phi_i\in w$ for every $1\leq i \leq n$, $\square\theta_i\in w$ for every $1\leq i\leq m$, and \begin{equation}\label{sh1} \vdash([\alpha]\phi_1\wedge\dots\wedge [\alpha]\phi_n)\wedge (\square\theta_1\wedge\dots\wedge\square\theta_m)\land \lnot[\alpha] \phi\to\bot. 
\end{equation} Let $\theta=\theta_1\land\dots \land\theta_m$. Since $\square$ distributes over conjunction, we have that $\square\theta \in w$. Now, \eqref{sh1} implies that \begin{equation}\label{sh2} \vdash([\alpha]\phi_1\wedge\dots\wedge [\alpha]\phi_n)\wedge \square\theta\to [\alpha] \phi.
\end{equation}

By necessitation of the ought-to-do operator, its axiom $(K)$, and its distributivity over conjunction, we get that \begin{equation}\label{ay2}\scriptsize\vdash(\odot [\alpha]([\alpha]\phi_1)\land\dots \land \odot [\alpha]([\alpha]\phi_n))\land \odot [\alpha](\square\theta) \to \odot [\alpha]([\alpha]\phi).\end{equation} Similarly as above, we have that $\odot [\alpha]([\alpha]\phi_i)\in w$ for every $1\leq i\leq n$. Therefore, their conjunction also lies in $w$. At the same time, $\square\theta\in w$ implies that $\square \square \theta\in w$ (by axiom $(4)$ for $\square$), and therefore that $\odot [\alpha] (\square \theta)\in w$. Therefore, \eqref{ay2} entails that $\odot [\alpha]([\alpha]\phi)\in w$. But axiom $(T)$ for $[\alpha]$, necessitation for the ought-to-do operator, and $(K)$ for the ought-to-do operator then imply that $\odot [\alpha]\phi\in w$, and this contradicts the previously shown fact that $\odot [\alpha]\phi\notin w$. 
\end{itemize}

Let $u$ be the $\Lambda$-MCS that includes $\Sigma_\alpha^{w}\cup \{\square\psi ; \square\psi \in w\} \cup \lnot [\alpha] \phi$. It is clear that $u\in \overline{w}$ (by the first item of \textbf{Claim 1}), that $\Sigma_\alpha^{w}\subseteq u$, and that $[\alpha]\phi\notin u$.
\item $(\Rightarrow)$ Assume that $\Sigma_\alpha^{w}\subseteq w$, and let $[\alpha ]\phi\in \Sigma_\alpha^{w}$. By assumption, we have that $[\alpha ]\phi\in w$. Since $w$ is a $\Lambda$-MCS, axiom $(4)$ for $[\alpha ]$ implies that $[\alpha ]([\alpha ]\phi)\in w$ as well. With the fact that $wR_\alpha^{\overline{w}} v$, this implies that $[\alpha ]\phi\in v$. Therefore $\Sigma_\alpha^{w} \subseteq v$ for every $v\in \mathtt{Choice}_\alpha(w)$.

$(\Leftarrow)$ Analogous.

\item \textbf{Claim 2.} For $w\in W^\Lambda$, the set $\{u\in \overline{w} ; \mathtt{Value}_\mathcal{O}(u)=1\}$ is not empty. 

(Proof of claim) 

Let $w\in W^\Lambda$. We will first show that $\bigcup_{\alpha\in Ags} \Sigma_\alpha^{w}$ is consistent. Let $w\in W^\Lambda$. Suppose that $\bigcup_{\alpha\in Ags} \Sigma_\alpha^{w}$ is not consistent. Then there must exist a set  $\{\phi_1,\dots,\phi_n\}$ of formulas of $\mathcal{L}_{\textsf{KO}}$ such that $\odot [\alpha_{i} ]\phi_i\in w$ for every $1\leq i \leq n$ and \begin{equation}\label{wayne} 
\vdash([\alpha_{1} ]\phi_1\wedge\dots\wedge [\alpha_{n}  ]\phi_n)\to \bot.
\end{equation}

Without loss of generality, we assume that $\alpha_i\neq\alpha_j$ for all $j\neq i$ such that $j,i\in \{1,\dots,n\}$ --this assumption hinges on the fact that any stit operator distributes over conjunction. Since $w$ is a $\Lambda$-MCS, then the fact that $\odot [\alpha_{i} ]\phi_i\in w$ for every $1\leq i \leq n$ implies that 
$(\odot [\alpha_{1} ]\phi_1\wedge\dots\wedge\odot [\alpha_{n} ]\phi_n)\in w$. Notice then that axiom $(Oic)$ yields that \begin{equation}\label{wayne1} 
\vdash(\odot [\alpha_{1} ]\phi_1\wedge\dots\wedge\odot [\alpha_{n} ]\phi_n)\to \Diamond [\alpha_{1} ]\phi_1\wedge\dots\wedge\Diamond[\alpha_{n} ]\phi_n.
\end{equation}

By axiom $(IA)$, we have that \begin{equation}\label{wayne2} 
\vdash\Diamond [\alpha_{1} ]\phi_1\wedge\dots\wedge\Diamond[\alpha_{n} ]\phi_n\to \Diamond( [\alpha_{1} ]\phi_1\wedge\dots\wedge[\alpha_{n} ]\phi_n).
\end{equation} Therefore, equations \eqref{wayne1}, \eqref{wayne2}, and \eqref{wayne}  imply that \begin{equation}\label{wayne3} 
\vdash(\odot [\alpha_{1} ]\phi_1\wedge\dots\wedge\odot [\alpha_{n} ]\phi_n)\to \Diamond \bot.
\end{equation}
But this is a contradiction, since we had that $(\odot [\alpha_{1} ]\phi_1\wedge\dots\wedge\odot [\alpha_{n} ]\phi_n)\in w$, and $w$ is a $\Lambda$-MCS. Therefore, $\bigcup_{\alpha\in Ags} \Sigma_\alpha^{w}$ is consistent.

Next, we show that the union $\bigcup_{\alpha\in Ags} \Sigma_\alpha^{w}\cup  \{\square\psi ; \square\psi \in w\}$ is also consistent. Suppose that this is not the case. Since $\bigcup_{\alpha\in Ags} \Sigma_\alpha^{w}$ is consistent, there must exist sets $\{\phi_1,\dots,\phi_n\}$ and $\{\theta_1,\dots,\theta_m\}$ of formulas of $\mathcal{L}_{\textsf{KO}}$ such that $\odot [\alpha_i ]\phi_i\in w$ for every $1\leq i \leq n$, $\square\theta_i\in w$ for every $1\leq i\leq m$, and \begin{equation}\label{sh'} \vdash([\alpha_1 ]\phi_1\wedge\dots\wedge [\alpha_n]\phi_n)\wedge (\square\theta_1\wedge\dots\wedge\square\theta_m)\to\bot. 
\end{equation} Let $\theta=\theta_1\land\dots \land\theta_m$. Since $\square$ distributes over conjunction, we have that  $\square\theta=\square\theta_1\wedge\dots\wedge\square\theta_m$.Since $w$ is a $\Lambda$-MCS, then $\square\theta\in w$. In these terms, \eqref{sh'} implies that \begin{equation}\label{sh''}\vdash([\alpha_1 ]\phi_1\wedge\dots\wedge [\alpha_n ]\phi_n)\to\lnot \square \theta.\end{equation} Once again, we assume without loss of generality that $\alpha_i\neq\alpha_j$ for all $j\neq i$ such that $j,i\in \{1,\dots,n\}$. Analogous to the procedure we used to show that $\bigcup_{\alpha\in Ags} \Sigma_\alpha^{w}$ is consistent, \eqref{sh''} implies that \begin{equation}\label{shorter} 
\vdash(\odot [\alpha_{1} ]\phi_1\wedge\dots\wedge\odot [\alpha_{n} ]\phi_n)\to \Diamond \lnot \square\theta.
\end{equation}

This entails that $\Diamond \lnot \square\theta\in w$, but this is a contradiction, since $\square\theta\in w$ implies with axiom $(4)$ for $square$ that $\square\square\theta\in w$. Let $u$ be the $\Lambda$-MCS that includes $\bigcup_{\alpha\in Ags} \Sigma_\alpha^{w}\cup  \{\square\psi ; \square\psi \in w\}$. It is clear that $u\in \overline{w}$ and that $\mathtt{Value}_\mathcal{O}(u)=1$. Therefore, the set $\{u\in \overline{w} ; \mathtt{Value}_\mathcal{O}(u)=1\}$ is not empty. Therefore, we have shown that \textbf{Claim 2.} is true. 

We proceed with the proof of the main statement. Recall that $\mathtt{Optimal}^{\overline{w}}_\alpha$ is defined as $\{L \in \mathtt{Choice}^{\overline{w}}_\alpha ; \textnormal{there is no } L' \in \mathtt{Choice}^{\overline{w}}_\alpha \textnormal{ such that }  L\prec L'\}$.

Let $w\in W^\Lambda$.
$(\Rightarrow)$ 
Let $w\in W^\Lambda$ such that $\Sigma_\alpha^{w}\subseteq w$. We will show that for every $L\in \mathtt{Choice}^{\overline{w}}_\alpha$ and $S\in \mathtt{State}^{\overline{w}}_\alpha$, $S\cap L \preceq S\cap \mathtt{Choice}^{\overline{w}}_\alpha(w)$, which would straightforwardly imply that $\mathtt{Choice}^{\overline{w}}_\alpha(w)\in \mathtt{Optimal}^{\overline{w}}_\alpha$. For this, it suffices to show that if $\mathtt{Value}_\mathcal{O}(u_*)=1$ for some $u_*\in S\cap L$, then $\mathtt{Value}_\mathcal{O}(o)=1$ for every $o\in S\cap \mathtt{Choice}^{\overline{w}}_\alpha(w)$. Therefore, take $L\in \mathtt{Choice}^{\overline{w}}_\alpha$, $S\in \mathtt{State}^{\overline{w}}_\alpha$, $o\in S\cap \mathtt{Choice}^{\overline{w}}_\alpha(w)$, and $u_*\in S\cap L$ such that $\mathtt{Value}_\mathcal{O}(u_*)=1$. This last thing means that $\bigcup_{\beta\in Ags} \Sigma_\beta^{w}\subseteq u_*$. 

As an intermediate result, we observe that for every $v\in S$, $\bigcup_{\beta\in Ags-\{\alpha\}} \Sigma_\beta^{w}\subseteq v$: let $v\in S$; since $u_*, v$ both lie in $S\in \mathtt{State}_\alpha^{\overline{w}}$, this means that $v\in \mathtt{Choice}^{\overline{w}}_\beta(u_*)$ for every $\beta\in Ags-\{\alpha\}$, so that item (b) above implies that $\Sigma_\beta^{w}\subseteq v$ for every $\beta\in Ags-\{\alpha\}$; therefore, $\bigcup_{\beta\in Ags-\{\alpha\}} \Sigma_\beta^{w}\subseteq v$. 

We want to show that $\mathtt{Value}_\mathcal{O}(o)=1$. Since $o$ was taken in $S\cap \mathtt{Choice}^{\overline{w}}_\alpha(w)$, then, on one hand, our intermediate result gives us that $\bigcup_{\beta\in Ags-\{\alpha\}} \Sigma_\beta^{w}\subseteq o$; on the other, we have that $o\in \mathtt{Choice}^{\overline{w}}_\alpha(w)$ implies by item (b) above that $\Sigma_\alpha^{w}\subseteq o$. With these two facts, we conclude that 
$\bigcup_{\beta\in Ags} \Sigma_\beta^{w}\subseteq o$, so that $\mathtt{Value}_\mathcal{O}(o)=1$. 

$(\Leftarrow)$  We work by contraposition. Suppose that $\Sigma_\alpha^{w}\not\subseteq w$. By item (b) above, this implies that $\Sigma_\alpha^{w}\not\subseteq o$ for every $o\in \mathtt{Choice}_\alpha^{\overline{w}} (w)$, so that $\mathtt{Value}_\mathcal{O}(o)=0$ for every $o\in \mathtt{Choice}^{\overline{w}}_\alpha(w)$. We will show that there exists an action in $\mathtt{Choice}^{\overline{w}}_\alpha$ that dominates $\mathtt{Choice}^{\overline{w}}_\alpha(w)$.

From \textbf{Claim 2.} we know that the set $\{u\in \overline{w} ; \mathtt{Value}_\mathcal{O}(u)=1\}$ is not empty. Take $u_*\in \{u\in \overline{u} ; \mathtt{Value}_{\mathcal{O}}(u)=1\}$. We will show that $\mathtt{Choice}_\alpha^{\overline{w}}(u_*)$ is the action that we are looking for, i.e., we will show ($\star$) that for every state $S\in \mathtt{State}^{\overline{w}}_\alpha$, $S\cap \mathtt{Choice}^{\overline{w}}_\alpha(w)\leq S\cap \mathtt{Choice}^{\overline{w}}_\alpha(u_*)$, and ($\star\star$) that there exists a state $S_*$ such that $S_*\cap \mathtt{Choice}^{\overline{w}}_\alpha(u_*)\nleq_s S_*\cap \mathtt{Choice}^{\overline{w}}_\alpha(w)$. For ($\star$), take any $S\in \mathtt{State}^{\overline{w}}_\alpha$. We want to show that $S\cap \mathtt{Choice}^{\overline{w}}_\alpha(w)\leq S\cap \mathtt{Choice}^{\overline{w}}_\alpha(u_*)$. Since for any $o\in S\cap \mathtt{Choice}^{\overline{w}}_\alpha(w)$ we have that $\mathtt{Value}_{\mathcal{O}}(o)=0$, this entails rather straightforwardly that $S\cap \mathtt{Choice}^{\overline{w}}_\alpha(w)\leq S\cap \mathtt{Choice}^{\overline{w}}_\alpha(u_*)$. For ($\star\star$), we set $S_{u_*}=\bigcap_{\beta\in Ags-\{\alpha\}}\mathtt{Choice}^{\overline{w}}_\beta(u_*)$. It is clear that $S_{u_*}\in \mathtt{State}^{\overline{w}}_\alpha$ and that $u_*\in S_{u_*}\cap \mathtt{Choice}^{\overline{w}}_\alpha(u_*)$. Similarly, we have that $(\mathtt{IA})_K$ implies that $S_{u_*}\cap \mathtt{Choice}^{\overline{w}}_\alpha(w)\neq \emptyset$. Since for any $o\in S_{u_*}\cap \mathtt{Choice}^{\overline{w}}_\alpha(w)$ we have that $\mathtt{Value}_{\mathcal{O}}(o)=0$, it is clear that $S_{u_*}\cap \mathtt{Choice}^{\overline{w}}_\alpha(u_*)\nleq S_{u_*}\cap \mathtt{Choice}^{\overline{w}}_\alpha(w)$. Now, our results ($\star$) and ($\star\star$) together render that  $\mathtt{Choice}^{\overline{w}}_\alpha(w)\prec \mathtt{Choice}^{\overline{w}}_\alpha(u_*)$, and this implies that $\mathtt{Choice}^{\overline{w}}_\alpha(w)\notin \mathtt{Optimal}^{\overline{w}}_\alpha$.

\item Let $w\in W^\Lambda$ and $L\in \mathtt{Choice}^{\overline{w}}_\alpha-\mathtt{Optimal}^{\overline{w}}_\alpha$. By \textbf{Claim 2.}, there exists $u_*\in \overline{w}$ such that $\mathtt{Value}_{\mathcal{O}}(u_*)=1$. Firstly, from item (c), we get that $\mathtt{Choice}^{\overline{w}}_\alpha(u_*)\in \mathtt{Optimal}^{\overline{w}}_\alpha$. The fact that $L\in \mathtt{Choice}^{\overline{w}}_\alpha-\mathtt{Optimal}^{\overline{w}}_\alpha$ implies with item (b) that $\mathtt{Value}_{\mathcal{O}}(w)=0$ for every $w\in L$, so that $S\cap L\leq S\cap \mathtt{Choice}^{\overline{w}}_\alpha(u_*)$ for every $S\in \mathtt{State}^{\overline{w}}_\alpha$ ($\star$). On the other hand, if we take $S_*$ to be the unique state in $\mathtt{State}^{\overline{w}}_\alpha$ such that $u_*\in S_*$, then it is clear that $S_*\cap \mathtt{Choice}^{\overline{w}}_\alpha(u_*)\nleq S_*\cap L$ ($\star\star$). Our results ($\star$) and ($\star\star$) together render that  $L\prec\mathtt{Choice}^{\overline{w}}_\alpha(u_*)$.

\end{enumerate}
\end{proof}

As for the previous results for $\odot_{\mathcal{S}}$, we first observe that for $\alpha\in Ags$, the $\mathbf{S5}$ axioms for the combination of operators $K_\alpha\square$ are theorems of $\Lambda$, which entails that the relation $R_\square \circ \approx_\alpha$ is an equivalence relation. For each $w\in W^\Lambda$, we will denote by $\alpha_k[w]$ the set $\{v\in W^\Lambda ; K_\alpha\square \phi\in w \Rightarrow \phi\in v\}$. Since the proof of the following lemma is the most important new development of the present work, we include its basic steps (the details are provided in the supplement).

\begin{lemma}[Existence for subjective ought] \label{kid} Fix $\alpha\in Ags$ and $w\in W^\Lambda$. \begin{enumerate}[(a)]
\item  For every formula $\phi$ of $\mathcal{L}_{\textsf{KO}}$, $\odot_\mathcal{S}[\alpha ]\phi\in w$ iff $K_\alpha\phi\in v$ for every $v\in \alpha_k[w]$ such that $\Gamma^{v}_\alpha\subseteq v$.
\item For every $v\in \sembrack{\mathtt{Choice}^{\overline{w}}_\alpha(w)}^{\overline{w'}}_\alpha$ ($w\approx_\alpha w'$), $\Gamma_\alpha^{w}\subseteq w$ iff $\Gamma_\alpha^{v} \subseteq v$. 
\item $\Gamma_\alpha^{w}\subseteq w$ iff $\mathtt{Choice}^{\overline{w}}_\alpha(w)\in \mathtt{S-Optimal}^{\overline{w}}_\alpha$. 
\item For each $L\in \mathtt{Choice}^{\overline{w}}_\alpha-(\mathtt{S-Optimal}^{\overline{w}}_\alpha)$, there exists $L'\in \mathtt{S-Optimal}^{\overline{w}}_\alpha$ such that $L\prec_s L'$. Thus, $\mathcal{M} ,w\Vdash \odot_\mathcal{S}[\alpha ]\phi  \textnormal{ iff } \forall L \in \mathtt{\mathtt{S-Optimal}}^{\overline{w}}_\alpha,  \forall w' \mbox{ such that } w\approx w', \ \mathcal{M} ,v\Vdash \phi \mbox{ for every } v\in \sembrack{L}^{\overline{w'}}_\alpha.$
\end{enumerate}
\end{lemma}
\begin{proof}
\begin{enumerate}[(a)]
\item $(\Rightarrow)$ We assume that $\odot_\mathcal{S}[\alpha ]\phi\in w$. Since $w$ is a $\Lambda$-MCS, this implies with axiom $(Cl)$ that $K_\alpha\square\odot_\mathcal{S}[\alpha ]\phi\in w$ \textcolor{blue}{}, so that the fact that $v\in \alpha_k[w]$ yields that $\odot_\mathcal{S}[\alpha ]\phi\in v$. The assumption that $\Gamma^{v}_\alpha\subseteq v$ then entails that $K_\alpha\phi\in v$. $(\Leftarrow)$ We work by contraposition. Suppose that $\odot_\mathcal{S}[\alpha ]\phi\notin w$. This implies that $\odot_\mathcal{S}[\alpha ]\phi\notin v$ for every $v\in \alpha_k[w]$, for if there existed $v\in \alpha_k[w]$ such that $\odot_\mathcal{S}[\alpha ]\phi\in v$, then $K_\alpha\square\odot_\mathcal{S}[\alpha ]\phi\in v$ by axiom $(Cl)$ and thus $\odot_\mathcal{S}[\alpha ]\phi\in w$, contradicting our assumption. Therefore, we take $v_*\in \alpha_k[w]$, where we know that $\odot_\mathcal{S}[\alpha ]\phi\notin v_*$. The union $\Gamma_\alpha^{v_*} \cup \{K_\alpha\square\psi; K_\alpha\square\psi \in w\}\cup \{\lnot K_\alpha \phi\}$ is consistent. Let $u$ be the $\Lambda$-MCS that includes $\Gamma_\alpha^{v_*} \cup \{K_\alpha\square\psi; K_\alpha\square\psi \in w\}\cup \{\lnot K_\alpha \phi\}$. We observe that for every $v\in \alpha_k[w]$, $\Gamma^{w}_\alpha=\Gamma^{v}_\alpha$. Then we have that $u\in \alpha_k[w]$, that $\Gamma_\alpha^u=\Gamma_\alpha ^{v_*}\subseteq u$, and that $K_\alpha\phi\notin u$.
\item Let $w'\in W^\Lambda$ such that $w\approx_\alpha w'$, and let $v\in \sembrack{\mathtt{Choice}^{\overline{w}}_\alpha(w)}^{\overline{w'}}_\alpha$. $(\Rightarrow)$  Assume that $\Gamma_\alpha^{w}\subseteq w$. Taking $v$ in $\sembrack{\mathtt{Choice}^{\overline{w}}_\alpha(w)}_\alpha^{\overline{w'}}$ implies that $w\approx_\alpha v$. Let $K_\alpha\phi\in \Gamma_\alpha^{w}$. By assumption, we have that $K_\alpha\phi\in w$. Axiom $(4)$ for $K_\alpha$ then implies that $K_\alpha K_\alpha\phi\in w$. Since $w\approx_\alpha v$, $K_\alpha\phi\in v$ as well. Therefore, we have shown that $\Gamma_\alpha^{w} \subseteq v$. $(\Leftarrow)$ Analogous.
\item $(\Rightarrow)$ Assume that $\Gamma_\alpha^{w}\subseteq w$. We show that for every $L\in \mathtt{Choice}^{\overline{w}}_\alpha$, $L\preceq_s \mathtt{Choice}^{\overline{w}}_\alpha(w)$. This comes from the fact that for any $w'$ such that such that $w\approx_\alpha w'$, if $\mathtt{Value}_\mathcal{S}(u_*)=1$ for some $u_*\in \sembrack{L}^{\overline{w'}}_\alpha\cap S$, then $\mathtt{Value}_\mathcal{S}(o)=1$ for every $o\in \sembrack{\mathtt{Choice}^{\overline{w}}_\alpha(w)}^{\overline{w'}}_\alpha\cap S$. $(\Leftarrow)$ We work by contraposition. Assume that $\Gamma_\alpha^{w}\not\subseteq w$.  The set $\{u\in \overline{w} ;\mathtt{Value}_{\mathcal{S}}(u)=1\}$ is not empty. Take $u_*\in \{u\in \overline{w} ;\mathtt{Value}_{\mathcal{S}}(u)=1\}$. It is the case that $\mathtt{Choice}^{\overline{w}}_\alpha(u_*)$ dominates $\mathtt{Choice}^{\overline{w}}_\alpha(w)$ in the order $\prec_{\mathcal{S}}$. 
\item Let $L\in \mathtt{Choice}^{\overline{w}}_\alpha-(\mathtt{S-Optimal}^{\overline{w}}_\alpha)$. We know that there exists $u_*\in \overline{w}$ such that $\mathtt{Value}_{\mathcal{S}}(u_*)=1$. From item (c) above, we get that $\mathtt{Choice}^{\overline{w}}_\alpha(u_*)\in \mathtt{S-Optimal}^{\overline{w}}_\alpha$. Notice that the fact that $L\in \mathtt{Choice}^{\overline{w}}_\alpha-(\mathtt{S-Optimal}^{\overline{w}}_\alpha)$ implies with items (b) and (c) above that that for every $w'$ such that $w\approx_\alpha w'$, $\mathtt{Value}_{\mathcal{S}}(o)=0$ for every $o\in \sembrack{L}^{\overline{w'}}_\alpha$. Therefore, it is clear that for every $w'$ such that $w\approx_\alpha w'$ and every $S\in \mathtt{State}^{\overline{w'}}_\alpha$, $\sembrack{L}^{\overline{w'}}_\alpha\cap S\leq_s \sembrack{\mathtt{Choice}^{\overline{w}}_\alpha(u_*)}^{\overline{w'}}_\alpha\cap S$ ($\star$). On the other hand, if we take $S_*$ to be the unique state in $\mathtt{State}^{\overline{w}}_\alpha$ such that $u_*\in S_*$, then $\sembrack{\mathtt{Choice}^{\overline{w}}_\alpha(u_*)}^{\overline{w}}_\alpha\cap S_*\nleq_s \sembrack{L}^{\overline{w}}_\alpha\cap S_*$ ($\star\star$). Our results ($\star$) and ($\star\star$) together render that  $L\prec_s\mathtt{Choice}^{\overline{w}}_\alpha(u_*)$.
\end{enumerate}
with this we conclude the proof.
\end{proof}

\begin{lemma}[Truth Lemma] \label{scof}
Let $\mathcal{M}$ be the canonical Kripke-estit model for $\Lambda$. For every formula $\phi$ of $\mathcal{L}_{\textsf{KO}}$ and every $w\in W^\Lambda$, we have that $\mathcal{M},w\Vdash \phi \ \textnormal{iff} \ \phi\in w.$
\end{lemma}

\begin{proof}
We proceed by induction on $\phi$. The cases with the boolean connectives  are standard. It remains to deal with the modal operators. For formulas involving $\square, [\alpha ]$, and $K_\alpha$ both directions follow straightforwardly from Lemma \ref{carajo} (1), (2), and (3), respectively. As for the operators for the objective and subjective oughts, we have the following arguments:
\begin{itemize}
\item (``$\odot[\alpha]$'')

($\Rightarrow$) We work by contraposition. We assume that $\odot[\alpha ]\phi\notin w$. By Lemma \ref{vang} (a), we get that there exists $v'\in \overline{w}$ such that $\Sigma_\alpha^{w}\subseteq v'$ and $[\alpha ]\phi\notin v'$. By Lemma \ref{carajo} (2), this implies that there exists $v\in \overline{w}$ such that $v'R^{\overline{w}}_\alpha v$ and such that $\phi\notin v$. By induction hypothesis, we get that $\mathcal{M} ,v\nVdash \phi$. From Lemma \ref{vang} (c), we know that $\mathtt{Choice}^{\overline{w}}_\alpha(v)=\mathtt{Choice}^{\overline{w}}_\alpha(v')\in \mathtt{Optimal}^{\overline{w}}_\alpha$, but this means that there exists an ``action'' in $\mathtt{Optimal}^{\overline{w}}_\alpha$ --namely $\mathtt{Choice}^{\overline{w}}_\alpha(v')$-- such that $v\in \mathtt{Choice}^{\overline{w}}_\alpha(v')$ and $\mathcal{M} ,v\nVdash \phi$, which implies by Lemma \ref{vang} (d) that $\mathcal{M} ,w\nVdash \odot [\alpha ]\phi$.

($\Leftarrow$) We assume that $\odot[\alpha ]\phi\in w$. By Lemma \ref{vang} (a), we get that $[\alpha ]\phi\in v'$ for every $v'\in \overline{w}$ such that $\Sigma_\alpha^{w}\subseteq v'$. Take $L\in \mathtt{Optimal}^{\overline{w}}_\alpha$. From Lemma \ref{vang} (b) we know that $\Sigma_\alpha^{w}\subseteq v$ for every $v\in L$. Therefore, for every $v\in L$, we have that $[\alpha ]\phi\in v$, which implies by axiom $(T)$ for the stit operator that for every $v\in L$, $\phi\in v$. By induction hypothesis, we get that for every $v\in L$, $\mathcal{M},v\Vdash\phi$, so that $\mathcal{M},w\Vdash\odot[\alpha ]\phi$.

\item (``$\odot_{\mathcal{S}}[\alpha]$'')

($\Rightarrow$) We work by contraposition. We assume that $\odot_{\mathcal{S}}[\alpha ]\phi\notin w$. By Lemma \ref{kid} (a), we get that there exists $v'\in \alpha_k[w]$ such that $\Gamma_\alpha^{v'}\subseteq v'$ and $K_\alpha \phi\notin v'$. By Lemma \ref{carajo} (3), this in turn implies that there exists $v\in W^\Lambda$ such that $v'\approx_\alpha v$ and such that $\phi\notin v$. Let $o_v\in \overline{w}$ be such that $v\approx_\alpha o_v$ (we know that such $o_v$ exists because condition $(\mathtt{US})_K$ guarantees that the fact that $v'\in\alpha_k[w]$ implies that there exists $o_v$ in $\overline{w}$ such that $v'\approx_\alpha o_v$, so that transitivity of $\approx_\alpha$ renders that  $v\approx_\alpha o_v$ as well). We observe that the fact that $o_v\approx_\alpha v$ implies that $v\in \sembrack{\mathtt{Choice}^{\overline{w}}_\alpha(o_v)}^{\overline{v}}_\alpha$. Now, by induction hypothesis, we get that  $\mathcal{M},v\nVdash \phi$. From Lemma \ref{kid} (b) we get that the facts that $\Gamma_\alpha^{v'}\subseteq v'$ and $v'\approx_\alpha o_v$ imply that $\Gamma_\alpha^{o_v}\subseteq o_v$, so that Lemma \ref{kid} (c) then implies $\mathtt{Choice}^{\overline{w}}_\alpha(o_v)\in \mathtt{S-Optimal}^{\overline{w}}_\alpha$. This means that there exists an ``action'' in $\mathtt{S-Optimal}^{\overline{w}}_\alpha$ --namely $\mathtt{Choice}^{\overline{w}}_\alpha(o_v)$-- such that $v\in \sembrack{\mathtt{Choice}^{\overline{w}}_\alpha(o_v)}^{\overline{v}}_\alpha$ and $\mathcal{M} ,v\nVdash \phi$, which implies by Lemma \ref{kid} (d) that $\mathcal{M},w\nVdash \odot_\mathcal{S}[\alpha ]\phi$.

($\Leftarrow$) We assume that $\odot_\mathcal{S}[\alpha ]\phi\in w$. By Lemma \ref{kid} (a), we get that $K_\alpha\phi\in v$ for every $v\in \alpha_k[w]$ such that $\Gamma_\alpha^{v}\subseteq v$. Take $L\in \mathtt{S-Optimal}^{\overline{w}}_\alpha$. From Lemma \ref{kid} (b) we know that for every $w'$ such that $w\approx_\beta w'$ , $\Gamma^{v'}_\alpha\subseteq v'$ for every $v'\in \sembrack{L}^{\overline{w'}}_\alpha$. Therefore, for every $v'\in \sembrack{L}^{\overline{w'}}_\beta$, we have that $K_\alpha \phi\in v'$, which implies --by axiom $(T)$ for the knowledge operator-- that for every $v'\in \sembrack{L}^{\overline{w'}}_\alpha$, $\phi\in v'$. By induction hypothesis, we get that for every $w'$ such that $w\approx_\beta w'$ and  $v'\in \sembrack{L}^{\overline{w'}}_\alpha$, $\mathcal{M},v'\Vdash\phi$. Since $L$ was taken as an arbitrary ``action'' in $\mathtt{S-Optimal}^{\overline{w}}_\alpha$, this means that $\mathcal{M} ,w\Vdash\odot_\mathcal{S}[\alpha ]\phi$, according to Lemma \ref{kid} (d).
\end{itemize}
With this we conclude the proof.
\end{proof}

\begin{theorem}[Completeness w.r.t. Kripke-estit models] \label{complete} The proof system $\Lambda$ is complete with respect to the class of Kripke-estit models.
\end{theorem}

\begin{proof}
Let $\phi$ be a $\Lambda$-consistent formula of $\mathcal{L}_{\textsf{KO}}$. Let $\Phi$ be the $\Lambda$-MCS including $\phi$. We have seen that the canonical Kripke-estit model is such that $\mathcal{M} , \Phi \Vdash \phi$.
\end{proof}

\begin{theorem} [Completeness w.r.t. bi-valued BT-models] \label{completeness} The proof system $\Lambda$ is complete with respect to the class of epistemic utilitarian bi-valued BT-models.
\end{theorem}

\begin{proof}
Let $\phi$ be a $\Lambda$-consistent formula of $\mathcal{L}_{\textsf{KO}}$. By Theorem \ref{complete} we know that there exists a Kripke-estit model $\mathcal{M}$ and a world $w$ in its domain such that $\mathcal{M} , w \Vdash \phi$. Proposition \ref{avion1} then ensures that the epistemic utilitarian bi-valued BT-model $\mathcal{M}^T$ associated to $\mathcal{M}$ is such that $\mathcal{M}^T, \langle\overline{w},h_w  \rangle\vDash \phi$. 
\end{proof}

\end{appendix}

\end{document}